\documentclass[10pt]{amsart}
\usepackage{amssymb, amsthm, amsmath}
\usepackage[all]{xy}
\usepackage{graphicx}
\usepackage{psfrag}

\theoremstyle{definition}

\newtheorem{example}{Example}

\newcommand\A{{\mathbb A}}
\newcommand{\B}{\mathbb B}
\newcommand{\D}{\mathbb D}
\newcommand{\E}{\mathbb E}

\newcommand\C{{\mathbb C}}
\newcommand\Q{{\mathbb Q}}
\newcommand\N{{\mathbb N}}

\newcommand{\cN}{{\mathcal N}}

\newcommand\IH{{\mathbb H}}
\newcommand{\ti}{\vartheta}

\newcommand{\Ti}{\Theta}

\newcommand{\om}{{\varpi}}
\newcommand\cW{{\mathcal W}}

\newcommand\cA{{\mathcal A}}

\newcommand\X{{\mathfrak X}}
\newcommand\x{{\mathrm{x}}}
\newcommand\Y{{\mathfrak Y}}
\newcommand\bb{{\omega}}

\newcommand\Z{{\mathbb Z}}

\newcommand\AS{{\mathfrak S}}

\newcommand\CS{{\mathfrak C}}
\newcommand\DS{{\mathfrak D}}

\newcommand\IGam{{\mathrm I}\Gamma}
\newcommand\ILam{{\mathrm I}\Lambda}
\newcommand\IB{{\mathrm I}B}

\newcommand\al{\alpha}
\newcommand{\be}{\beta}
\newcommand\la{\lambda}

\newcommand\Eta{H}

\newcommand\ssm{\smallsetminus}
\newcommand\noin{\noindent}

\newcommand\bull{{\scriptscriptstyle \bullet}}
\newcommand\eqto{\stackrel{\lower1.5pt\hbox{$\scriptstyle\sim\,$}}\to}
\newcommand\ov{\overline}

\newcommand\wh{\widehat}
\newcommand\wt{\widetilde}

\DeclareMathOperator{\Pf}{Pfaffian}
\DeclareMathOperator{\Sp}{Sp}
\DeclareMathOperator{\SO}{SO}
\DeclareMathOperator{\SL}{SL}
\DeclareMathOperator{\GL}{GL}
\DeclareMathOperator{\LG}{LG}
\DeclareMathOperator{\IG}{IG}

\DeclareMathOperator{\OG}{OG}
\DeclareMathOperator{\G}{G}

\DeclareMathOperator{\HH}{\mathrm{H}}

\DeclareMathOperator{\type}{\mathrm{type}}

\newcommand{\ignore}[1]{}
\newcommand{\pic}[2]{\includegraphics[scale=#1]{#2}}

\begin{document}

\title[Theta and eta polynomials]
{Theta and eta polynomials in geometry, Lie theory, and combinatorics}

\date{April 14, 2020}

\author{Harry~Tamvakis} 
\address{University of Maryland, Department of
Mathematics, William E. Kirwan Hall, 4176 Campus Drive, 
College Park, MD 20742, USA}
\email{harryt@math.umd.edu}



\begin{abstract}
The classical Schur polynomials form a natural basis for the ring of
symmetric polynomials and have geometric significance since they
represent the Schubert classes in the cohomology ring of
Grassmannians.  Moreover, these polynomials enjoy rich combinatorial
properties. In the last decade, an exact analogue of this picture has
emerged in the symplectic and orthogonal Lie types, with the Schur
polynomials replaced by the theta and eta polynomials of Buch, Kresch,
and the author. This expository paper gives an overview of what is
known to date about this correspondence, with examples.
\end{abstract}

\maketitle

\section{Introduction}

The structure of the cohomology ring of the Grassmannian $\G(m,n)$ of
$m$ dimensional linear subspaces of $\C^n$ was first explored by
Schubert \cite{Sc}, Pieri \cite{Pi}, and Giambelli \cite{G}.  The ring
$\HH^*(\G(m,n),\Z)$ has an additive basis of Schubert classes, coming
from the cell decomposition of $\G(m,n)$, and there is a natural
choice of multiplicative generators for this ring, namely, the special
Schubert classes. Giambelli showed that a general Schubert class,
which may be indexed by a partition, can be expressed as a
Jacobi-Trudi determinant \cite{J, Tr} whose entries are special
classes. It soon became apparent (see for example \cite{Le}) that the
resulting algebra is closely connected to the theory of Schur
polynomials. The latter polynomials were originally defined by Cauchy
\cite{C}, and studied by many others since then, motivated to a large
extent by their applications to the representation theory of the
symmetric and general linear groups \cite{S1}, and related
combinatorics.

The theory of theta and eta polynomials, by contrast, has its origins
in geometry, and specifically in the desire to extend the
aforementioned work of Giambelli to the cohomology of symplectic and
orthogonal Grassmannians. The first steps in this direction were taken
in the 1980s by Hiller and Boe \cite{HB} and Pragacz \cite{P}. They
proved Pieri and Giambelli formulas for the Grassmannians of maximal
isotropic subspaces, with the Schubert classes indexed by strict
partitions and the Jacobi-Trudi determinants replaced by Schur
Pfaffians \cite{S2}.

In 2008, Buch, Kresch, and the author announced a series of works
\cite{BKT1, BKT2, BKT3, T3} which went beyond these hermitian
symmetric (or cominuscule) examples. Two crucial insights from
op.\ cit.\ were the identification of the correct set of special
Schubert class generators to employ, and the realization of the
essential role that Young's {\em raising operators} \cite{Y} play in
the theory. The latter becomes clear only when one attempts to
understand the cohomology of {\em nonmaximal} isotropic
Grassmannians.  Our papers introduced {\em $k$-strict partitions} to
index the Schubert classes and {\em theta polynomials} to represent
them, in both the classical and quantum cohomology rings of general
symplectic and odd orthogonal Grassmannians. Three years later, the
companion papers \cite{BKT4, T6} dealt with the even orthogonal case, using
{\em typed $k$-strict partitions} and {\em eta polynomials}.

Theta and eta polynomials can be viewed as symmetric polynomials for
the action of the corresponding Weyl group, but this is not obvious
from their definition, and was pointed out only recently \cite{T9}. In
fact, a substantial part of the theory of Schur polynomials can be
extended to the world of theta and eta polynomials, but this requires
a change in perspective, as well as the introduction of new techniques
of proof. It turns out that these objects can be applied to solve the
(equivariant) Giambelli problem for the classical Lie groups, that is,
to obtain intrinsic polynomial representatives for the (equivariant)
Schubert classes on any classical $G/P$-space -- so, any (isotropic)
partial flag manifold \cite{BKTY, T5, T7}. The resulting formulas
(\ref{Apartfl}), (\ref{Cpartfl}), and (\ref{Dpartfl}) are stated
using solely the language of Lie theory.

The goal of this expository paper is to illustrate the correspondence
between Schur polynomials and theta/eta polynomials in the case of
single polynomials, where the story is most complete. There is ample
room for further interesting connections to be found, and the reader
is likely to discover more by just asking for the theta/eta
analogue of their favorite statement about Schur polynomials.  We have
stopped short of discussing extensions of some of these results to the
theory of degeneracy loci and equivariant cohomology, quantum
cohomology, and $K$-theory. 

This article is organized as follows. The Schur, theta, and eta
polynomials are featured in Sections \ref{schurps}, \ref{thetapolys},
and \ref{etapolys}, respectively. Each of these is split into parallel
subsections on initial definitions and Pieri rules, the cohomology of
Grassmannians, symmetric polynomials, algebraic combinatorics, and the
cohomology of partial flag manifolds. Finally, Section \ref{notes}
contains historical notes and references.

It is a pleasure to thank the Hellenic Mathematical Society for the
invitation to lecture in the First Congress of Greek Mathematicians in
Athens, June 2018. The present paper is a more detailed and expanded
version of my talk there.

\section{Schur polynomials}
\label{schurps}

\subsection{Definition using raising operators}
\label{first}
We will define the classical Schur polynomials by using the
Jacobi-Trudi formula, but rewritten in the language of Young's raising
operators. Let $u_1, u_2,\ldots$ be a sequence of commuting
independent variables, and set $u_0:=1$ and $u_i:=0$ for
$i<0$. Throughout the paper, $\al:=(\al_1,\al_2,\ldots)$ will denote
an integer sequence with only finitely many nonzero terms
$\al_i$. For any such $\al$, we let $u_\al:=u_{\al_1}u_{\al_2}\cdots$,
which is a monomial in the variables $u_i$.

Given any integer sequence $\alpha$ and $i<j$, we define
\[
R_{ij}(\alpha) := (\alpha_1,\ldots,\alpha_i+1,\ldots,\alpha_j-1,
\ldots).
\] 
A {\em raising operator} $R$ is any monomial in the $R_{ij}$'s. For 
any such $R$, we let $R\,u_{\al} := u_{R(\al)}$, where, as is customary, 
we regard the operator $R$ as acting on the index $\al$, and not on 
the monomial $u_\al$. 

An integer sequence $\al$ is a {\em partition} if $\al_i\geq
\al_{i+1}\geq 0$ for each $i\geq 1$. If $\la$ is a partition, we let
$|\la|$ denote the sum of all its parts $\la_i$, and the length
$\ell(\la)$ be the number of $i$ such that $\la_i\neq 0$.  A partition
$\la$ may be represented using its Young diagram of boxes, placed in
left justified rows, with $\la_i$ boxes in the $i$-th row, for each
$i\geq 1$. The diagram of the conjugate partition $\wt{\la}$ is
obtained by taking the transpose of the diagram of $\la$.  The
containment relation $\la\subset\mu$ between two partitions is defined
using their respective diagrams; in this case the set-theoretic
difference $\mu\ssm\la$ is the skew diagram $\mu/\la$. For example,
we illustrate below the diagram of the partition $(4,3,3,2)$ and the
skew diagram $(4,3,3,2)/(3,1,1)$.
\medskip
\[
\includegraphics[scale=0.28]{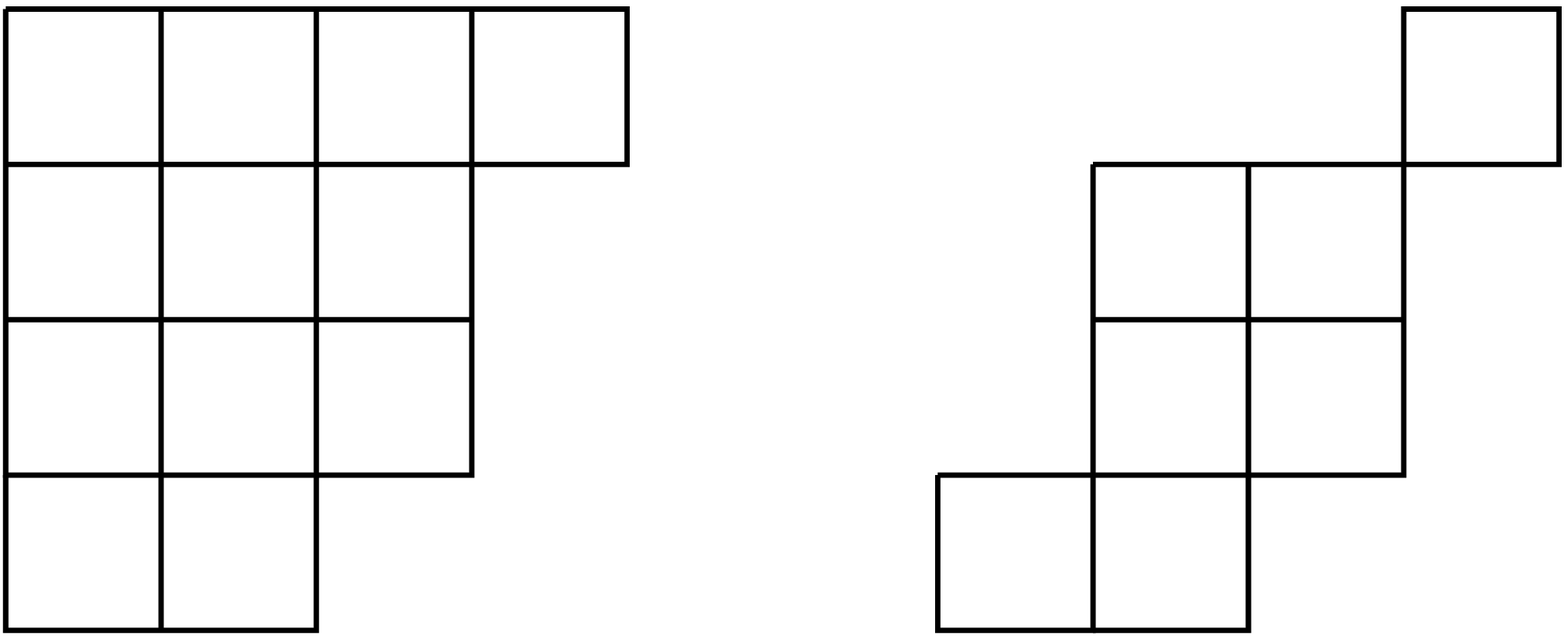}
\]
\smallskip

The {\em Schur polynomial} $s_\la(u)$ is 
defined by the raising operator formula
\begin{equation}
\label{giambelliA}
s_{\la}(u) := \prod_{i<j} (1-R_{ij})\, u_{\la}.
\end{equation}
We regard the products in formulas such as (\ref{giambelliA}) as 
being over all pairs $(i,j)$ with $j\leq N$ for some $N$ greater than 
or equal to the length $\ell(\la)$. The result will always be independent
of any such $N$, so we may assume that $N=\ell(\la)$.

\begin{example}
\label{ex1}
For any partition $\la=(a,b)$ with two parts $a$ and $b$, we have
\[
s_{a,b}(u)=(1-R_{12}) \, u_{a,b} = 
u_{a,b} - u_{a+1,b-1} = u_au_b-u_{a+1}u_{b-1} = 
\left|\begin{array}{cc}
u_a & u_{a+1} \\ u_{b-1} & u_b
\end{array}\right|.
\]
\end{example}

The formula of Example \ref{ex1} generalizes: for any partition $\la$
of length $\ell$, we have
\begin{equation}
\label{equality1}
s_\la(u) =  \det(u_{\la_i+j-i})_{1\leq i,j\leq \ell}.
\end{equation}
The equivalence of (\ref{giambelliA}) with (\ref{equality1}) is a formal
consequence of the Vandermonde identity
\[
\prod_{1\leq i<j\leq \ell}(\x_i-\x_j)= \det(\x_i^{\ell-j})_{1\leq i,j\leq \ell} 
\]

The monomials $u_\la$ and the polynomials $s_\la(u)$ as $\la$ runs
over all partitions form two $\Z$-bases of the graded polynomial ring
$\A:=\Z[u_1,u_2,\ldots]$. In addition to the Giambelli formula
(\ref{giambelliA}), these two bases of $\A$ interact via the {\em
  Pieri rule}. To state the latter, recall that a skew diagram is a
{\em horizontal strip} (resp.\ a {\em vertical strip}) if it does not
contain two boxes in the same column (resp.\ row).  For any partition
$\la$ and $p\geq 0$, we then have the Pieri rule
\[
u_p\cdot s_\la(u) = \sum_{\mu} s_\mu(u)
\]
with the sum over all partitions $\mu\supset \la$ such that $\mu/\la$
is a horizontal strip with $p$ boxes.

\begin{example} We have
\[
u_3\cdot s_{2,2,1}(u) = s_{5,2,1}(u)+s_{4,2,2}(u)+s_{4,2,1,1}(u)+
s_{3,2,2,1}(u).
\]
\end{example}

\subsection{Cohomology of Grassmannians} 

Let $\G(m,n)$ be the Grassmannian of all $m$-dimensional complex
linear subspaces of $\C^n$. The general linear group $\GL_n(\C)$ acts
transitively on $\G(m,n)$, and the stabilizer of the point $\langle
e_1,\ldots, e_m\rangle$ under this action is the parabolic subgroup
$P_m$ of matrices in $\GL_n(\C)$ in the block form
\[
\left(
\begin{array}{c|c}
* & * \\  \hline
0 & *
\end{array}
\right)
\]
where the lower left block is an $(n-m)\times m$ zero matrix. 
It follows that 
\[
\G(m,n)=\GL_n/P_m
\]
and hence that $\G(m,n)$ is a complex manifold of dimension $m(n-m)$. 
Furthermore, $\G(m,n)$ has the structure of an algebraic variety,
and the same is true of all the geometric objects which appear in this paper.

The Grassmannian $\G(m,n)$ has a natural decomposition into Schubert
cells $\X^\circ_\la$, one for each partition $\la$ whose diagram is
contained in an $m\times (n-m)$ rectangle. The Schubert variety
$\X_\la$ is the closure of the cell $\X^\circ_\la$, and is an algebraic
subvariety of $\G(m,n)$ of complex codimension equal to $|\la|$.
Concretely, if $e_1,\ldots,e_n$ is the canonical basis of 
$\C^n$, and $F_i$ is the $\C$-linear span of $e_1,\ldots,e_i$, then
\[
\X_\la := \{V\in \G(m,n)\ |\ \dim(V\cap F_{n-m+i-\la_i})\geq i, \ 1\leq i \leq m\}.
\]
If $[\X_\la]$ denotes the cohomology class of $\X_\la$, then the cell
decomposition of $\G(m,n)$ implies that there is an isomorphism of
abelian groups
\begin{equation}
\label{Grasseq}
\HH^*(\G(m,n),\Z) \cong \bigoplus_\la \Z [\X_\la].
\end{equation}

For every integer $p$ with $1\leq p \leq n-m$, the variety $\X_p$ is
the locus of subspaces $V$ in $\G(m,n)$ which meet the subspace
$F_{n-m+1-p}$ nontrivially.  The varieties $\X_p$ are the {\em
  special Schubert varieties}, and their cohomology classes $[\X_p]$
are the {\em special Schubert classes}.  Let $Q\to \G(m,n)$ denote the
universal quotient vector bundle over $\G(m,n)$. Then for each
integer $p\geq 0$, the $p$-th Chern class $c_p(Q)$ of $Q$ is equal to
$[\X_p]$ in $\HH^{2p}(\G(m,n),\Z)$.  We can now state the {\em
  Giambelli formula}
\begin{equation}
\label{GiambA}
[\X_\la] = s_{\wt{\la}}(c(Q))
\end{equation}
where $\wt{\la}$ is the conjugate partition of $\la$, and the Chern
class polynomial $s_{\wt{\la}}(c(Q))$ is obtained from
$s_{\wt{\la}}(u)$ by performing the substitutions $u_p\mapsto c_p(Q)$
for each integer $p$.  The discussion in \S \ref{first} and
equation (\ref{GiambA}) imply that the {\em Pieri rule}
\begin{equation}
\label{PieriA}
[\X_p]\cdot [\X_\la] = \sum_{\mu} [\X_\mu]
\end{equation}
also holds in $\HH^*(\G(m,n),\Z)$, where the sum is over all indexing
partitions $\mu$ containing $\la$ such that $\mu/\la$ is a horizontal
strip with $p$ boxes.

\subsection{Symmetric polynomials}
\label{spsA}

Fix an integer $n\geq 1$ and let $X_n:=(x_1,\ldots,x_n)$, where the 
$x_i$ are independent variables. The Weyl group $S_n$ of $\GL_n$ 
acts on the polynomial ring $\Z[X_n]$ by permuting the variables, 
and the invariant subring is the ring $\Lambda_n:=\Z[X_n]^{S_n}$
of symmetric polynomials. Two important families of elements of 
$\Lambda_n$ are the elementary symmetric polynomials $e_p(X_n)$ and 
the complete symmetric polynomials $h_p(X_n)$. These are defined 
the by the generating function equations 
\[
\sum_{p=0}^\infty e_p(X_n)t^p = \prod_{i=1}^n(1+x_it) \quad \mathrm{and}
\quad 
\sum_{p=0}^\infty h_p(X_n)t^p = \prod_{i=1}^n(1-x_it)^{-1},
\]
respectively, where $t$ is a formal variable. The fundamental theorem of 
symmetric polynomials states that 
\[
\Lambda_n=\Z[e_1(X_n),\ldots,e_n(X_n)].
\]

For each partition $\la$, the Schur polynomial $s_\la(X_n)$ is obtained
from $s_\la(u)$ by making the substitution $u_p\mapsto h_p(X_n)$ for every
integer $p$. In this way, we obtain the {\em Jacobi-Trudi formula}
\[
s_\la(X_n)=\prod_{i<j} (1-R_{ij})\, h_\la(X_n) = 
\det(h_{\la_i+j-i}(X_n))_{1\leq i,j\leq \ell(\la)},
\]
where we have set $h_\al:=h_{\al_1}h_{\al_2}\cdots$ 
for every integer sequence $\al$. If we similarly let 
$e_\al:=e_{\al_1}e_{\al_2}\cdots$, then we have the dual Jacobi-Trudi formula
\begin{equation}
\label{JTdual}
s_\la(X_n)=\prod_{i<j} (1-R_{ij})\, e_{\wt{\la}}(X_n) = 
\det(e_{\wt{\la}_i+j-i}(X_n))_{1\leq i,j\leq \ell(\wt{\la})}.
\end{equation}
We deduce that 
\[
\Lambda_n=\bigoplus_\la \Z \,s_\la(X_n)
\]
where the sum is over all partitions $\la=(\la_1,\ldots,\la_n)$ of 
length at most $n$.

The classical definition of Schur polynomials is as a quotient of
alternant determinants. For each integer $k\geq 1$, let
$\delta_k:=(k,\ldots,1,0)$, and, for any integer vector
$\al=(\al_1,\ldots,\al_n)$, let $x^{\al}:=x_1^{\al_1}\cdots
x_n^{\al_n}$.  Define the alternating operator $\cA$ on $\Z[X_n]$ by
\[
\cA(f):=\sum_{\om\in S_n} (-1)^{\ell(\om)}\om(f),
\]
where $\ell(\om)$ denotes the length of the permutation $\om$. 
Then we have
\begin{equation}
\label{Cauchydef}
s_\la(X_n)=\left. \det(x_i^{\la_i+n-j})_{i,j}\right\slash \det(x_i^{n-j})_{i,j} = 
\left. \cA(x^{\la+\delta_{n-1}})\right\slash \cA(x^{\delta_{n-1}}).
\end{equation}

We restate equation (\ref{Cauchydef}) in a form closer to its analogue
for theta polynomials. For each $r\geq 1$, embed $S_r$ into $S_{r+1}$
by adjoining the fixed point $r+1$. Let $S_\infty=\cup_rS_r$ denote
the corresponding infinite symmetric group of bijections $\om:\Z_{>0}
\to \Z_{>0}$ such that $\om_i=i$ for all but finitely many $i$. Here,
and in the sequel, $\om_i$ denotes the value $\om(i)$, for each $i
\geq 1$.  The {\em code} of $\om$ is the the sequence
$\gamma=\gamma(\om)$ with $\gamma_i:=\#\{j>i\ |\ \om_j<\om_i\}$. The
     {\em shape} of $\om$ is the partition $\la=\la(\om)$ whose parts
     are the nonzero entries $\gamma_i$ arranged in weakly decreasing
     order.

\begin{example}
An $n$-Grassmannian permutation $\om$ is an element of $S_\infty$ such
that $\om_i < \om_{i+1}$ for each $i\neq n$. The shape of any such
$\om$ is the partition $\la:=(\om_n-n,\ldots,\om_1-1)$. Conversely, any
partition $\la=(\la_1,\ldots,\la_n)$ of length at most $n$ corresponds
to a unique $n$-Grassmannian permutation $\om$ with $\la(\om)=\la$.
\end{example}

Let $\om$ be an $n$-Grassmannian element of $S_\infty$ with
corresponding partition $\la(\om)$. If $\om_0:=(n,n-1,\ldots,1)$
denotes the longest permutation in $S_n$, then observe that
$\la(\om_0)=\delta_{n-1}$ and $\la(\om\om_0)=\la(\om)+\delta_{n-1}$.
Therefore, we have
\[
s_{\la(\om)}(X_n)=\left. \cA(x^{\la(\om\om_0)})\right\slash \cA(x^{\la(\om_0)}).
\]

\subsection{Algebraic combinatorics}
\label{ac}

Although not strictly necessary, in this subsection we will extend our
family of variables to be an infinite sequence $X:=(x_1,x_2,\ldots)$, 
and work with formal power series instead of polynomials. Define
$h_p(X)$ using the generating function expansion
\[
\sum_{p=0}^{\infty}h_p(X)t^p = \prod_{i=1}^\infty (1-x_it)^{-1},
\]
and, for any partition $\la$, the Schur function
\[
s_\la(X):=\prod_{i<j}(1-R_{ij})\, h_\la(X).
\]

A {\em tableau} $T$ on the shape $\la$
is a filling of the boxes of $\la$ with positive integers, so that
the entries are weakly increasing along each row from left to right
and strictly increasing down each column. The content vector 
$c(T)=(n_1,n_2,\ldots)$ of $T$ has $n_j$ equal to the number of 
entries $j$ in $T$. We then have the {\em tableau formula}
\begin{equation}
\label{tabform}
s_\la(X)= \sum_Tx^{c(T)}
\end{equation}
with the sum over all tableaux $T$ of shape $\la$. Equation (\ref{tabform})
shows that $s_\la(X)$ is a formal power series with nonnegative integer
coefficients, which have a combinatorial interpretation.

\begin{example}
A Young tableaux of shape $\la=(2,1)$ is of the form
\[
\includegraphics[scale=0.23]{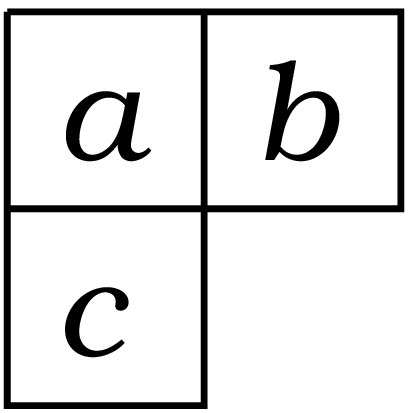}
\]
where the positive integers $a,b,c$ satisfy $a\leq b$ and $a<c$. We
therefore have
\[
s_{2,1}(X)= \sum_{{a\leq b} \atop {a<c}}x_ax_bx_c = 
\sum_{a\neq b}x_a^2x_b + 2\sum_{a<b<c}x_ax_bx_c.
\]
\end{example}

In order to represent the Schubert classes not only on Grassmannians, but 
on any partial flag variety, we require a generalization of (\ref{tabform}) 
which involves certain symmetric functions defined by Stanley. 
The group $S_\infty$ is generated by the simple transpositions 
$s_i=(i,i+1)$ for $i\geq 1$. A {\em reduced word} of a permutation $\om$
in $S_\infty$ is a sequence $a_1\cdots a_\ell$ of positive integers
such that $\om=s_{a_1}\cdots s_{a_\ell}$ and $\ell$ is minimal, so (by
definition) equal to the length $\ell(\om)$ of $\om$. 

The {\em nilCoxeter algebra} $\cN_n$ of the symmetric group $S_n$ is
the free associative algebra with unit generated by the elements
$\xi_1,\ldots,\xi_{n-1}$, modulo the relations
\[
\begin{array}{rclr}
\xi_i^2 & = & 0 & i\geq 1\ ; \\
\xi_i\xi_j & = & \xi_j\xi_i & |i-j|\geq 2\ ; \\
\xi_i\xi_{i+1}\xi_i & = & \xi_{i+1}\xi_i\xi_{i+1} & i\geq 1.
\end{array}
\]
For any $\om\in S_n$, choose a reduced word $a_1\cdots a_\ell$ for
$\om$ and define $\xi_\om := \xi_{a_1}\ldots \xi_{a_\ell}$. Then the $\xi_\om$
for $\om\in S_n$ are well defined, independent of the choice of reduced 
word $a_1\cdots a_\ell$, and form a free $\Z$-basis of
$\cN_n$. We denote the coefficient of $\xi_\om\in \cN_n$ in the
expansion of the element $\zeta\in \cN_n$ by $\langle \zeta,\om\rangle$.
We therefore have 
\[
\zeta = \sum_{\om\in S_n}\langle \zeta,\om\rangle\,\xi_\om,
\]
for all $\zeta\in \cN_n$.

Recall that $t$ denotes an indeterminate and define
\begin{gather*}
A(t) := (1+t \xi_{n-1})(1+t \xi_{n-2})\cdots 
(1+t \xi_1)
 \ ; \\
A(X) := A(x_1)A(x_2)\cdots
\end{gather*}
and a function $G_\om(X)$ for $\om\in S_n$ by
\[
G_\om(X) := \langle A(X), \om\rangle. 
\]
It turns out that $G_\om(X)$ is symmetric in the $x_i$ variables; it 
is called a {\em Stanley symmetric function}. Clearly, $G_\om$
has nonnegative integer coefficients. 

When $\om$ is the $n$-Grassmannian permutation associated to a
partition $\la$ of length at most $n$, then $G_\om(X)=s_\la(X)$. More
generally, for any permutation $\om$, when $G_\om$ is expanded in the
basis of Schur functions, we have
\begin{equation}
\label{Geq}
G_\om(X) = \sum_{\la\, :\, |\la| = \ell(\om)}c^\om_\la \,s_\la(X)
\end{equation}
for some nonnegative integers $c^\om_\la$.  That is, the symmetric
function $G_\om$ is {\em Schur positive}. There exist several
different combinatorial interpretations of the coefficients
$c^\om_\la$, but the most important one for our purposes uses the {\em
transition trees} of Lascoux and Sch\"utzenberger, as subsequently defined.

We say that a permutation $\om=(\om_1,\om_2,\ldots)$ has a {\em
  descent} at position $i\geq 1$ if $\om_i>\om_{i+1}$. For $i<j$, let
$t_{ij}$ denote the transposition which interchanges $i$ and $j$.  For
any permutation $\om\in S_\infty$, the {\em transition tree} $T(\om)$
of $\om$ is constructed as follows.  The tree $T(\om)$ is a rooted
tree with nodes given by permutations of the same length $\ell(\om)$,
and root $\om$.  If $\om=1$ or $\om$ is Grassmannian, then set
$T(\om):=\{\om\}$. Otherwise, let $r$ be the largest descent of $\om$,
and set 
\[
s := \max(j>r\ |\ \om_j < \om_r).
\] 
Define
\[
I(\om):=\{i \ |\ 1\leq i < r \ \ \mathrm{and} \ \ \ell(\om
t_{rs}t_{ir}) = \ell(\om) \}
\]
and let
\[
\Psi(\om):=\begin{cases}
\{ \om t_{rs}t_{ir}\ |\ i\in I(\om)\} & \text{if $I(\om)\neq \emptyset$}, \\
\Psi(1\times \om) & \text{otherwise}.
\end{cases}
\]
We define $T(\om)$ recursively, by
joining $\om$ by an edge to each $v\in \Psi(\om)$, and attaching to
each $v\in \Psi(\om)$ its tree $T(v)$. One can show that $T(\om)$ is a
finite tree whose leaves are all Grassmannian permutations. Moreover,
the Stanley coefficient $c^\om_\la$ in equation (\ref{Geq}) is equal
to the number of leaves of shape $\la$ in the transition tree $T(\om)$
associated to $\om$.

\begin{example}
The transition tree for the permutation $\om=(2,1,5,4,3)$ is shown below.
\medskip
\[
\includegraphics[scale=0.31]{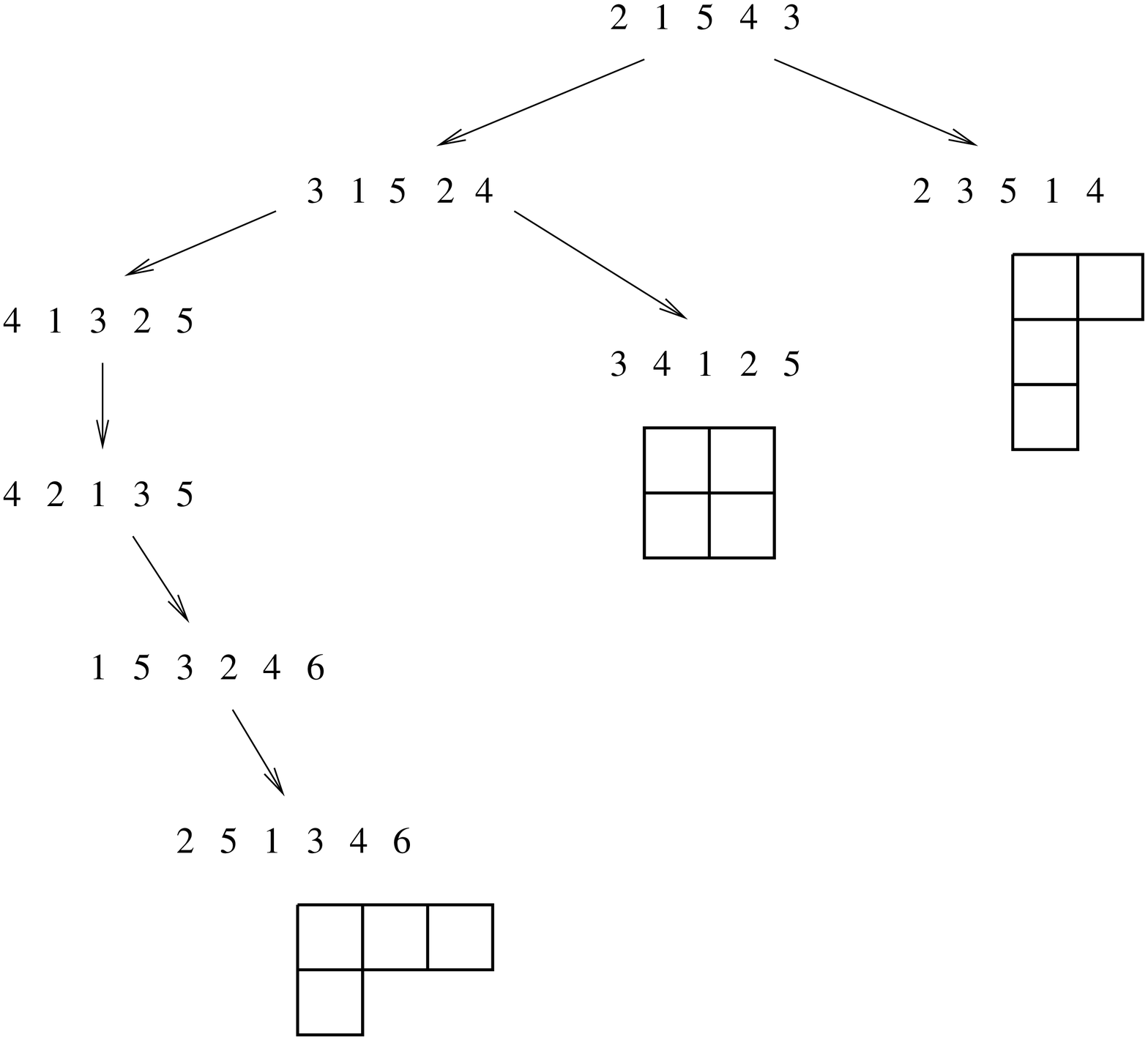}
\]
\medskip
\medskip
It follows that
\[
G_{21543}(X)=s_{3,1}(X)+s_{2.2}(X)+s_{2,1,1}(X).
\]
\end{example}

\subsection{Cohomology of flag manifolds}

Let $\{e_1,\ldots,e_n\}$ denote the standard basis of
$E:=\C^n$ and let $F_i=\langle e_1,\ldots, e_i\rangle$ 
be the subspace spanned by the first $i$ vectors of this basis.
The group $G=\GL_n(\C)$ acts transitively on the space of all complete 
flags in $E$, and the stabilizer of the flag $F_\bull$ is the Borel
subgroup $B$ of invertible upper triangular matrices. If $T\subset B$
denotes the maximal torus of diagonal matrices, then the Weyl group 
$W=S_n$ can be identified with $N_G(T)/T$.

The parabolic subgroups $P$ of $\GL_n$ with $P\supset B$
correspond to sequences $a_1<\cdots < a_p$ of positive
integers with $a_p<n$. For any such $P$, the manifold $\GL_n/P$
parametrizes partial flags of subspaces
\[
0 = E_0 \subset E_1 \subset \cdots \subset E_p \subset E=\C^n
\]
with $\dim(E_r) = a_r$ for each $r\in [1,p]$. We agree that
$E_r$ and $E$ will also denote the corresponding
tautological vector bundles over $\GL_n/P$. The associated
parabolic subgroup $W_P$ of $S_n$ is generated by the
simple transpositions $s_i$ for $i\notin\{a_1,\ldots, a_p\}$. 

There is a canonical presentation of the cohomology ring of $\GL_n/B$,
which gives geometric significance to the variables which appear in
Section \ref{spsA}. Let $\ILam_n$ denote the ideal of $\Z[X_n]$
generated by the homogeneous elements of positive degree in
$\Lambda_n$, so that $\ILam_n=\langle
e_1(X_n),\ldots,e_n(X_n)\rangle$. We then have ring isomorphism
\begin{equation}
\label{BorelA}
\HH^*(\GL_n/B) \cong \Z[X_n]/\ILam_n
\end{equation}
which maps each variable $x_i$ to $-c_1(E_i/E_{i-1})$. Moreover, for any 
parabolic subgroup $P$ of $\GL_n$, the projection map $\GL_n/B\to
\GL_n/P$ induces an injection $\HH^*(\GL_n/P)\hookrightarrow
\HH^*(\GL_n/B)$ of cohomology rings, and we have
\begin{equation}
\label{BorelPA}
\HH^*(\GL_n/P) \cong \Z[X_n]^{W_P}/\ILam_n^P,
\end{equation}
where $\Z[X_n]^{W_P}$ denotes the $W_P$-invariant subring of $\Z[X_n]$, 
and $\ILam_n^P$ is the ideal of $\Z[X_n]^{W_P}$ generated by 
$e_1(X_n),\ldots,e_n(X_n)$.

Consider the set
\[
W^P := \{\om\in S_n\ |\ \ell(\om s_i) = \ell(\om)+1,\  \forall\, i \notin
\{a_1,\ldots, a_p\}, \ i<n\}
\]
of minimal length $W_P$-coset representatives in $S_n$. We have
a decomposition
\[
\GL_n = \bigcup_{\om \in W^P}B\om P
\]
and for each $\om \in W^P$, the $B$-orbit of $\om P$ in $\GL_n/P$ is
the {\em Schubert cell} $\Y^\circ_\om:=B\om P/P$. Let $\Y_\om$ be the
closure of $\Y^\circ_\om$ in $\GL_n/P$, and set
$\X_\om:=\Y_{\om_0\om}$.  The {\em Schubert class} $[\X_\om]$ is the
cohomology class of $\X_\om$ in $\HH^{2\ell(\om)}(\GL_n/P,\Z)$. We
thus obtain an isomorphism of abelian groups
\[
\HH^*(\GL_n/P,\Z) \cong \bigoplus_{\om\in W^P}\Z[\X_\om]
\]
which generalizes (\ref{Grasseq}). 

If $V_1$ and $V_2$ are complex vector bundles over a manifold $M$, define
cohomology classes $\phi_p$ by the generating function equation
\[
\sum_{p=0}^\infty \phi_pt^p = c_t(V_2^*)/c_t(V_1^*), 
\]
where $c_t(V_i^*)=1-c_1(V_i)t+c_2(V_i)t^2-\cdots$ is the Chern
polynomial of $V_i^*$, for $i=1,2$. Given any partition $\la$, the
polynomial $s_\la(V_1-V_2)$ is obtained from $s_{\la}(u)$ via the
substitution $u_p\mapsto \phi_p$ for every $p\in\Z$, so that
\[
s_\la(V_1-V_2):= \prod_{i<j} (1-R_{ij})\, \phi_\la.
\]

Recall that $E_r$ for $r\in [1,p]$ and $E$ denote the tautological and
trivial rank $n$ vector bundles over $\GL_n/P$, respectively. For
any $\om\in W^P$, we then have
\begin{equation}
\label{Apartfl}
[\X_\om] = \sum_{\underline{\la}} c^\om_{\underline{\la}}\,
s_{\wt{\la}^1}(E-E_1) s_{\wt{\la}^2}(E_1-E_2)\cdots s_{\wt{\la}^p}(E_{p-1}-E_p)
\end{equation}
in $\HH^*(\GL_n/P,\Z)$, where the sum is over all sequences of
partitions $\underline{\la}=(\la^1,\ldots,\la^p)$ and the 
coefficients $c^\om_{\underline{\la}}$ are given by 
\begin{equation}
\label{Acoeff}
c^\om_{\underline{\la}} := \sum_{u_1\cdots u_p = \om}
c_{\la^1}^{u_1}\cdots c_{\la^p}^{u_p}
\end{equation}
summed over all factorizations $u_1\cdots u_p= \om$ such that
$\ell(u_1)+\cdots + \ell(u_p)=\ell(\om)$ and $u_j(i)=i$ for all $j>1$
and $i \leq a_{j-1}$. The nonnegative integers $c^{u_i}_{\la^i}$ which
appear in the summands in (\ref{Acoeff}) agree with the Stanley
coefficients from equation (\ref{Geq}).  When $p=1$, the partial flag
manifold $\GL_n/P$ is the Grassmannian $\G(a_1,n)$, and formula
(\ref{Apartfl}) specializes to equation (\ref{GiambA}).

\begin{example}
\label{ASex}
Let $P=B$ be the Borel subgroup, so that the flag manifold $\GL_n/B$
parametrizes complete flags of subspaces $0=E_0\subset E_1 \subset
\cdots \subset E_n =\C^n$. For each $i\in [1,n]$, let
$x_i:=-c_1(E_i/E_{i-1})$.  Then for any partition $\la$, we have
\[
s_{\wt{\la}}(E_{i-1}-E_i) = 
\begin{cases}
x_i^r & \text{if $\la=r\geq 0$}, \\
0 & \text{otherwise}.
\end{cases}
\]
Let $A_i(t):=(1+t \xi_{n-1})(1+t \xi_{n-2})\cdots (1+t \xi_i)$ in
$\cN_n[t]$, employing the notation of Section \ref{ac}. Define the {\em
  Schubert polynomial} $\AS_\om(X_n)$ by
\begin{equation}
\label{FSeq}
\AS_\om(X_n):= \langle A_1(x_1)\cdots A_{n-1}(x_{n-1}), \om\rangle.
\end{equation}
It is then straightforward to show that formula (\ref{Apartfl}) is 
equivalent to the statement that, for any permutation $\om\in S_n$, 
we have $[\X_\om]=\AS_\om(X_n)$ in $\HH^*(\GL_n/B,\Z)$.
\end{example}

\section{Theta polynomials} 
\label{thetapolys}

\subsection{Definition and Pieri rule}
\label{firstC}

Fix a nonnegative integer $k$. A partition $\la$ is called {\em
  $k$-strict} if no part greater than $k$ is repeated, that is,
$\la_i>k\Rightarrow \la_i >\la_{i+1}$. A {\em strict} partition is the same
as a $0$-strict partition. For a general $k$-strict partition $\la$,
we define the operator
\[
R^{\la} := \prod_{i<j} (1-R_{ij})\prod_{\la_i+\la_j > 2k+j-i}
(1+R_{ij})^{-1}
\]
where the first product is
over all pairs $i<j$ and second product is over pairs $i<j$ such that
$\la_i+\la_j > 2k+j-i$. The {\em theta polynomial} $\Ti^{(k)}_\la(u)$ of level 
$k$ is defined by 
\begin{equation}
\label{giambelliC}
\Ti^{(k)}_{\la}(u) := R^\la u_{\la}.
\end{equation}
We will write $\Ti_\la(u)$ for $\Ti^{(k)}_\la(u)$ when the level $k$ is
understood.

\begin{example}
\label{ex1C}
(a) Suppose that $\la=(a,b)$ has two parts $a$ and $b$ with $a+b>2k+1$.
Then we have
\begin{align*}
\Ti_{a,b}(u) &= \frac{1-R_{12}}{1+R_{12}} \, u_{a,b} = 
(1-2R_{12}+2R_{12}^2-2R_{12}^3+\cdots) \, u_{a,b} \\
&= u_au_b-2u_{a+1}u_{b-1} +2u_{a+2}u_{b-2}-2u_{a+3}u_{b-3}+\cdots.
\end{align*}

\medskip
\noin
(b) If $\la_i\leq k$ for each $i$, then 
\begin{equation}
\label{Schurdet}
\Theta_\la(u) = \prod_{i<j}(1-R_{ij})\, u_\la  = 
 \det(u_{\la_i+j-i})_{1\leq i,j\leq \ell(\la)}.
\end{equation}

\medskip
\noin
(c) If $\la_i>k$ for all nonzero parts $\la_i$, then
\begin{equation}
\label{SchurPf}
\Theta_\la(u) = \prod_{i<j}\frac{1-R_{ij}}{1+R_{ij}}\,u_{\la} = 
\Pf\left(\Ti_{\la_i,\la_j}(u)\right)_{1\leq i<j\leq \ell'}
\end{equation}
where $\ell'$ is the least positive even integer such that 
$\ell'\geq \ell(\la)$.
\end{example}

Example \ref{ex1C} shows that as $\la$ varies, the theta polynomial
$\Ti_\la$ interpolates between the determinant (\ref{Schurdet}) and
the Pfaffian (\ref{SchurPf}). The equality of Example \ref{ex1C}(c) is
a formal consequence of Schur's Pfaffian identity
\[
\prod_{1\leq i<j \leq \ell'}\frac{\x_i-\x_j}{\x_i+\x_j} =
\Pf\left(\frac{\x_i-\x_j}{\x_i+\x_j}\right)_{1\leq i,j \leq \ell'}.
\]

\begin{example}
Let $k=2$ and $\la=(5,2,1)$. Then we have
\begin{align*}
\Ti^{(2)}_{5,2,1}(u) & = 
\frac{1-R_{12}}{1+R_{12}}(1-R_{13})(1-R_{23}) \, u_{5,2,1} \\
&= (1-2R_{12}+2R_{12}^2 - 2 R_{12}^3)(1-R_{13}-R_{23}) 
\, u_{5,2,1} \\
&= u_5u_2u_1 -u_5u_3 -2 u_6u_1^2 + u_6u_2+2u_7u_1.
\end{align*}
\end{example}

Recall that $\A=\Z[u_1,u_2,\ldots]$, and let $\A^{(k)}$ be the quotient
of $\A$ by the ideal of relations
\begin{equation}
\label{relations}
\frac{1-R_{12}}{1+R_{12}}\, u_{p,p} = 
u_p^2 + 2\sum_{i=1}^p(-1)^i u_{p+i}u_{p-i}= 0
\ \ \ \text{for} \  p > k.
\end{equation}
Then the monomials $u_\la$ and the polynomials $\Ti_\la(u)$ as $\la$ runs
over all $k$-strict partitions form two $\Z$-bases of the graded ring
$\A^{(k)}$. To state the Pieri rule for the $\Ti_\la(u)$, which holds 
modulo the relations (\ref{relations}), we need some further definitions. 

We say that the box in row $r$ and column $c$ of a $k$-strict
partition $\lambda$ is {\em $k$-related} to the box in row $r'$ and
column $c'$ if $|c-k-1|+r = |c'-k-1|+r'$.  For example, the two grey
boxes in the following partition are $k$-related.
\[ \pic{0.65}{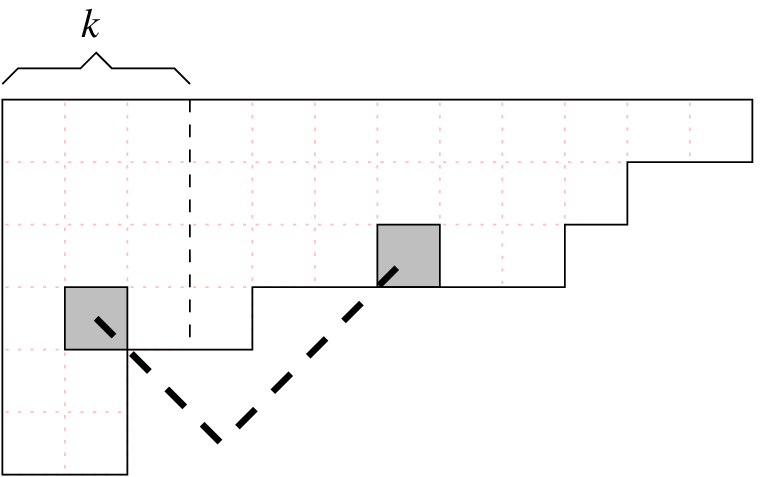} \]

For any two $k$-strict partitions $\lambda$ and $\mu$, we have a
relation $\lambda \xrightarrow{p} \mu$ if $|\mu|=|\la|+p$ and $\mu$ is
obtained by removing a vertical strip from the first $k$ columns of
$\lambda$ and adding a horizontal strip to the resulting diagram, so
that

(1) if one of the first $k$ columns of $\mu$ has the same number of
boxes as the same column of $\lambda$, then the bottom box of this
column is $k$-related to at most one box of $\mu \smallsetminus
\lambda$; and

(2) if a column of $\mu$ has fewer boxes than the same column of
$\lambda$, then the removed boxes and the bottom box of $\mu$ in this
column must each be $k$-related to exactly one box of $\mu
\smallsetminus \lambda$, and these boxes of $\mu \smallsetminus
\lambda$ must all lie in the same row.

If $\lambda \xrightarrow{p} \mu$, we let $\D$ be the set of boxes of
$\mu\ssm \la$ in columns $k+1$ and higher which are {\em not}
mentioned in (1) or (2). Define $N(\lambda,\mu)$ to be the number
of connected components of $\D$ which do not have a box in column
$k+1$. Here we consider that two boxes are connected if they have at 
least one vertex in common.

For any $k$-strict partition $\la$ and $p\geq 0$, we then
have the {\em Pieri rule}
\begin{equation}
\label{pruleC}
u_p\cdot \Ti_\la(u) = \sum_{\mu} 2^{N(\la,\mu)}\Ti_\mu(u)
\end{equation}
in $\A^{(k)}$, where the sum over all $k$-strict partitions $\mu$ such that
$\la\xrightarrow{p}\mu$.

\begin{example} 
\label{ex2C}
(a) When $k=1$, we have the following equality in $\A^{(1)}$:
\[
u_3\cdot \Ti_{2,1}(u) = 2\, \Ti_6(u) + 4\,\Ti_{5,1}(u) + \Ti_{4,2}(u)+
2\,\Ti_{4,1,1}(u) +\Ti_{3,2,1}(u).
\]
 
\noin
(b) If $|\la|+p \leq k$, then 
\begin{equation}
\label{Pieriksmall}
u_p\cdot \Ti_\la(u) = \sum_{\mu} \Ti_\mu(u)
\end{equation}
holds in $\A^{(k)}$, where the sum is over all partitions $\mu\supset
\la$ such that $|\mu|=|\la|+p$ and $\mu/\la$ is a horizontal strip.

\medskip
\noin
(c) If $\la_i>2k$ for all nonzero parts $\la_i$ and $p\geq 0$ is arbitrary, 
then
\begin{equation}
\label{Pieriklarge}
u_p\cdot \Ti_\la(u) = \sum_{\mu} 2^{N(\la,\mu)}\Ti_\mu(u)
\end{equation}
holds in $\A^{(k)}$, where the sum is over all strict partitions $\mu\supset
\la$ such that $|\mu|=|\la|+p$ and $\mu/\la$ is a horizontal strip,
and $N(\la,\mu)$ equals the number of connected components of
$\mu/\la$ which do not meet the first column.
\end{example}

In harmony with Example \ref{ex1C}, Example \ref{ex2C} illustrates
that as $p$ and $\la$ vary, the rule for the product $u_p\cdot
\Ti_\la$ interpolates between the Pieri rule (\ref{Pieriksmall}) for
Schur polynomials and (\ref{Pieriklarge}), which is the Pieri rule for
the Schur $Q$-functions (see Example \ref{exQ}).

\subsection{Cohomology of Grassmannians} 
\label{cGC}

Equip the vector space $\C^{2n}$ with the nondegenerate skew-symmetric
bilinear form $(\ ,\,)$ defined by the conditions $(e_i,e_j)=0$ for
$i+j\neq 2n+1$ and $(e_i,e_{2n+1-i})=1$ for $1\leq i \leq n$. The
symplectic group $\Sp_{2n}(\C)$ is the subgroup of $\GL_{2n}(\C)$
consisting of those elements $g$ such that $(gv_1,gv_2)=(v_1,v_2)$,
for every $v_1,v_2\in \C^{2n}$. We say that a linear subspace $V$ of
$\C^{2n}$ is {\em isotropic} if the restriction of $(\ ,\,)$ to $V$
vanishes identically. Since the form is nondegenerate, we have
$\dim(V)\leq n$ for any isotropic subspace $V$. If $V$ is isotropic
and $\dim(V)=n$ then we call $V$ a {\em Lagrangian} subspace.

Fix an integer $k$ with $0\leq k \leq n-1$. The isotropic Grassmannian
$\IG=\IG(n-k,2n)$ parametrizes all isotropic linear subspaces of
dimension $n-k$ in $\C^{2n}$. The group $\Sp_{2n}$ acts
transitively on $\IG(n-k,2n)$, and the stabilizer of a fixed isotropic
$(n-k)$-plane under this action is a maximal parabolic subgroup $P_k$
of $\Sp_{2n}$, so that $\IG(n-k,2n) = \Sp_{2n}/P_k$.

The Schubert cells $\X^\circ_\la$ in $\IG(n-k,2n)$ are indexed by
the $k$-strict partitions whose diagrams are contained in an
$(n-k)\times (n+k)$ rectangle. The Schubert variety $\X_\la$ is the
closure of the $\X^\circ_\la$, and has codimension $|\la|$ in $\IG$.
If $F_i$ denotes the $\C$-linear span of $e_1,\ldots,e_i$ for each
$i\in [1,2n]$, then
\[
\X_\lambda := \{V \in \IG \mid \dim(V \cap
   F_{p_j(\lambda)}) \geq j \ \ \ \forall\, 1 \leq j \leq n-k \},
\]
where the strictly increasing index function $\{p_j(\lambda)\}_{1\leq
  j \leq n-k}$ is defined by
\[
p_j(\lambda) := n+k+j-\lambda_j - \#\{i<j \ |\  \la_i+\la_j>2k+j-i \}.
\]
If $[\X_\la]$ denotes the cohomology class of $\X_\la$ in
$\HH^{2|\la|}(\IG,\Z)$, then we have a group isomorphism
\begin{equation}
\label{IGrasseq}
\HH^*(\IG(n-k,2n),\Z) \cong \bigoplus_\la \Z [\X_\la].
\end{equation}

The varieties $\X_p$ for $1\leq p \leq n+k$ are the {\em special
  Schubert varieties}.
Let $Q\to \IG(n-k,2n)$ denote the universal quotient vector bundle
over $\IG$, which has rank $n+k$. For every integer $p\geq 0$, the
$p$-th Chern class $c_p(Q)$ is equal to the {\em special Schubert
  class} $[\X_p]$ in $\HH^{2p}(\IG(n-k,2n),\Z)$. We now have the {\em
  Giambelli formula}
\begin{equation}
\label{GiambC}
[\X_\la] = \Ti_{\la}(c(Q))
\end{equation}
where the Chern class polynomial $\Ti_{\la}(c(Q))$ is obtained from
$\Ti_{\la}(u)$ by performing the substitutions $u_p\mapsto c_p(Q)$
for each integer $p$. Moreover, the {\em Pieri rule} 
\begin{equation}
\label{PieriC}
[\X_p]\cdot [\X_\la] = \sum_{\la\xrightarrow{p}\mu} 2^{N(\la,\mu)}[\X_\mu]
\end{equation}
holds in $\HH^*(\IG(n-k,2n),\Z)$, where the sum is over all partitions 
$\mu$ such that $\la\xrightarrow{p}\mu$ and the diagram of $\mu$ fits 
in an $(n-k)\times (n+k)$ rectangle.

The ring $\A^{(k)}$ is naturally isomorphic to the stable cohomology
ring
\[
\IH(\IG_k) = \lim_{\longleftarrow}\HH^*(\IG(n-k,2n),\Z)
\]
of the isotropic Grassmannian $\IG$. This is the inverse limit in
the category of {\em graded} rings of the directed system
\[
\cdots \leftarrow \HH^*(\IG(n-k,2n),\Z) \leftarrow
\HH^*(\IG(n+1-k,2n+2),\Z) \leftarrow \cdots
\]
Under this isomorphism, the variables $u_p$ map to the Chern classes
$c_p(Q)$ of the universal quotient bundle $Q$ over $\IG$. If $S$
denotes the tautological subbundle of the trivial rank $2n$ vector
bundle over $\IG(n-k,2n)$, then the symplectic form $(\ ,\,)$ gives a
pairing $S\otimes Q\to {\mathcal O}$, which in turn produces an
injection $S \hookrightarrow Q^*$. Using the Whitney sum formula
$c(S)c(Q)=1$, we therefore obtain
\[
c(Q^*)c(Q) = c(Q^*)c(S)^{-1}=c(Q^*/S).
\]
The relations (\ref{relations}) arise from this, 
using the fact that $c_r(Q^*/S)=0$ for $r>2k$.

\subsection{Symmetric polynomials}
\label{spsC}

Let $c:=(c_1,c_2,\ldots)$ be a sequence of commuting independent
variables, set $c_0:=1$ and $c_p=0$ for $p<0$, and for every integer
sequence $\al$, let $c_\al:=c_{\al_1}c_{\al_2}\cdots$. For any raising
operator $R$, we let $R\, c_\al:=c_{R\al}$. Consider the graded ring
$\Gamma$ which is the quotient of the polynomial ring $\Z[c]$ modulo
the ideal generated by the relations
\[
c_p^2+2\sum_{i=1}^p(-1)^ic_{p+i}c_{p-i}=0, \ \ \ \text{for all $p\geq 1$}.
\]
The ring $\Gamma$ is isomorphic to $\A^{(0)}$ and to the stable cohomology ring
\[
\IH(\LG) = \lim_{\longleftarrow}\HH^*(\LG(n,2n),\Z)
\]
of the Lagrangian Grassmannian $\LG$, with the variables
$c_p$ mapping to the Chern classes
$c_p(Q)$ of the universal quotient bundle $Q\to\LG$.

The Weyl group of the symplectic group $\Sp_{2n}$ is the {\em
  hyperoctahedral group} of signed permutations on the set
$\{1,\ldots,n\}$, which is the semidirect product $S_n \ltimes \Z_2^n$
of $S_n$ with $\Z_2^n=\{\pm 1\}^n$. We use a bar over an entry to
denote a negative sign; thus $w=(\ov{2},\ov{3},1)$ maps $(1,2,3)$ to
$(-2,-3,1)$. The group $W_n$ is generated by the simple transpositions
$s_i=(i,i+1)$ for $1\leq i \leq n-1$ and the sign change $s_0$ which
satisfies $s_0(1)=\ov{1}$ and $s_0(j)=j$ for all $j\geq 2$. A {\em
  reduced word} of an element $w\in W_n$ is a sequence $a_1\cdots
a_\ell$ of nonnegative integers of minimal length $\ell$ such that
$w=s_{a_1}\cdots s_{a_\ell}$. The number $\ell$ is called the {\em
  length} of $w$, and denoted by $\ell(w)$.

Fix $n\geq 1$ and let $X_n:=(x_1,\ldots,x_n)$, as in Section \ref{spsA}.
There is a natural action of $W_n$ on $\Gamma[X_n]$ which extends the 
action of $S_n$ on $\Z[X_n]$, defined as follows. 
The simple transpositions $s_i$ for $i\in [1,n-1]$ 
act by interchanging $x_i$ and $x_{i+1}$ while leaving all the remaining
variables fixed. The reflection $s_0$ satisfies $s_0(x_1)=-x_1$ and 
$s_0(x_j)=x_j$ for all $j \geq 2$, while
\begin{equation}
\label{equivaC}
s_0(c_p) := c_p+2\sum_{j=1}^p x_1^jc_{p-j} \ \ \text{for all $p\geq 1$}.
\end{equation}
If $t$ denotes a formal variable which is fixed by $s_0$, then we 
express equation (\ref{equivaC}) using generating functions as
\[
s_0\left(\sum_{p=0}^\infty c_pt^p\right) = \frac{1+x_1t}{1-x_1t}
\left(\sum_{p=0}^\infty c_pt^p\right).
\]

For every integer $p$, define an element ${}^nc_p$ of $\Gamma[X_n]$ by
\[
{}^nc_p:=\sum_{j=0}^p c_{p-j}e_j(X_n), \ \ \ \text{for $p\geq 1$}
\]
and let $\Gamma^{(n)} := \Z[{}^nc_1, {}^nc_2, \ldots]$. Let
$\Gamma[X_n]^{W_n}$ denote the subring of $W_n$-invariants in
$\Gamma[X_n]$. We claim that the generators ${}^nc_p$ of
$\Gamma^{(n)}$ lie in $\Gamma[X_n]^{W_n}$. Indeed, we clearly have
$s_j({}^nc_p)={}^nc_p$ for each $j\geq 1$, while
\begin{align*}
s_0\left(\sum_{p=0}^\infty {}^nc_pt^p\right) &=
s_0\left(\sum_{p=0}^\infty c_pt^p \cdot \prod_{j=1}^n(1+x_jt)\right)
\\ &= \frac{1+x_1t}{1-x_1t} \left(\sum_{p=0}^\infty c_pt^p
\right)\cdot (1-x_1t)\prod_{j=2}^n(1+x_jt) = \sum_{p=0}^\infty
      {}^nc_pt^p.
\end{align*}
In fact, there is an equality
\begin{equation}
\label{Cinv}
\Gamma[X_n]^{W_n}=\Gamma^{(n)}= \Z[{}^nc_1, {}^nc_2, \ldots].
\end{equation}

The map which sends $u_p$ to ${}^nc_p$ for every integer $p$ induces a ring 
isomorphism $\A^{(n)}\cong \Gamma^{(n)}$. We therefore have
\begin{equation}
\label{GTieq}
\Gamma^{(n)}=\bigoplus_\la \Z\, \Ti_\la(X_n)
\end{equation}
where the sum is over all $n$-strict partitions $\la$, and the
polynomial $\Ti_\la(X_n)$ is obtained from the theta polynomial
$\Ti^{(n)}_\la(u)$ by making the substitution $u_p\mapsto {}^nc_p$ for
all $p$. In other words, we have
\[
\Ti_\la(X_n)=\Ti^{(n)}_\la(X_n) := R^\la\, ({}^nc)_\la,
\]
where we set $({}^nc)_\al:={}^nc_{\al_1}{}^nc_{\al_2}\cdots$ for any
integer sequence $\al$, and the raising operators $R_{ij}$ in $R^\la$
act on the subscripts $\al$ as usual.

For each $r\geq 1$, we embed $W_r$ in $W_{r+1}$ by adding the element
$r+1$ which is fixed by $W_r$, and set $W_\infty :=\cup_r W_r$.  Let
$w\in W_\infty$ be a signed permutation.  To simpify the notation, 
let $w_i$ denote the value $w(i)$, for each $i\geq 1$. Define a strict
partition $\mu(w)$ whose parts are the absolute values of the negative
entries of $w$, arranged in decreasing order. Let the {\em A-code} of
$w$ be the sequence $\gamma$ with $\gamma_i:=\#\{j>i\ |\ w_j<w_i\}$,
and define a partition $\delta(w)$ whose parts are the nonzero
entries $\gamma_i$ arranged in weakly decreasing order. The {\em
  shape} of $w$ is the partition $\la(w):=\mu(w)+\nu(w)$, where
$\nu(w):=\wt{\delta(w)}$ is the conjugate of $\delta(w)$. One can
show that the length $\ell(w)$ of $w$ is equal to $|\la(w)|$.

\begin{example}
(a) An $n$-Grassmannian signed permutation $w$ is an element of $W_\infty$ such
that $w_1>0$ and $w_i < w_{i+1}$ for each $i\neq n$. The shape of any
such $w$ is the $n$-strict partition $\la$ satisfying
\[
\la_i=\begin{cases} 
n+|w_{n+i}| & \text{if $w_{n+i}<0$}, \\
\#\{r\leq n\, :\, w_r> w_{n+i}\} & \text{if $w_{n+i}>0$}.
\end{cases}
\]
Conversely, any $n$-strict partition $\la$ corresponds to a unique
$n$-Grassmannian permutation $w$ with $\la(w)=\la$. 

\medskip
\noin
(b) Let $w_0:=(\ov{1},\ldots,\ov{n})$ be
the longest element of $W_n$. Then $\mu(w_0)=\delta_n$, $\nu(w_0)=
\delta_{n-1}$, and $\la(w_0)=\delta_n+\delta_{n-1}=(2n-1,2n-3,\ldots,
1)$.
\end{example}

Let $w$ be an $n$-Grassmannian element of $W_\infty$ with
corresponding partition $\lambda(w)$, and let $w_0$ be the longest
element of $W_n$.  For any integer sequence $\al$ and composition
$\beta$, let ${}^\be c_\al:=
{}^{\be_1}c_{\al_1}{}^{\be_2}c_{\al_2}\cdots$, and set $R_{ij}{}^\be
c_\al:= {}^\be c_{R_{ij}\al}$ for each $i<j$. Consider the {\em
  multi-Schur Pfaffian}
\begin{equation}
\label{multiSPf}
{}^{\nu(ww_0)}Q_{\la(ww_0)} := \prod_{i<j} \frac{1-R_{ij}}{1+R_{ij}}
\,{}^{\nu(ww_0)}c_{\la(ww_0)}.
\end{equation}
Define the alternating operator $\cA'$ on $\Gamma[X_n]$ by
\[
\cA'(f):=\sum_{w\in W_n} (-1)^{\ell(w)}w(f),
\]
where $\ell(w)$ is the length of the signed permutation $w$. 
We then have
\begin{equation}
\label{mScheq}
\Ti_{\la(w)}(X_n) 
=(-1)^{n(n+1)/2}\left. 
\cA'\left({}^{\nu(ww_0)}Q_{\la(ww_0)}\right)\right\slash
\cA'\left(x^{\la(w_0)}\right)
\end{equation}
in $\Gamma[X_n]$.

\subsection{Algebraic combinatorics}
\label{acC}

The combinatorial formulas discussed in this section require another
incarnation of the ring $\Gamma$, using the formal power series known
as Schur $Q$-functions. Fix a nonnegative integer $k$, set
$X_k:=(x_1,\ldots,x_k)$ and let $Z:=(z_1,z_2,\ldots)$ be a sequence of
variables. For any integer $p$, define $\ti_p=\ti_p(Z\,;X_k)$ by the
generating function equation
\begin{equation}
\label{tieq}
\sum_{p=0}^\infty \ti_pt^p = \prod_{i=1}^\infty \frac{1+z_it}{1-z_it}
\prod_{j=1}^k(1+x_jt).
\end{equation}
By definition, for every $k$-strict partition $\la$,
the theta polynomial $\Ti_\la(Z\,;X_k)$ is obtained from
$\Ti^{(k)}_\la(u)$ by making the substitution $u_p\mapsto \ti_p$
for every integer $p$. In other words, we have
\begin{equation}
\label{Thetatheta}
\Ti_\la(Z\,;X_k):= R^\la\, \ti_\la
\end{equation}
where, for any integer sequence $\al$,
$\ti_\al=\ti_{\al_1}\ti_{\al_2}\cdots$ and the raising operators
$R_{ij}$ in $R^\la$ are applied to $\ti_\la$ as usual.  Note that
$\Ti_\la(Z\,;X_k)$ is a formal power series in the $Z$ variables and a
polynomial in the variables $x_1,\ldots,x_k$.

\begin{example}
\label{exQ}
Suppose that $k=0$ and $\la$ is a strict partition. Then the formal
power series $Q_\la(Z):=\Ti_\la(Z)$ is a {\em Schur $Q$-function}. The
map which sends $c_p$ to $Q_p:=Q_p(Z)$ for every integer $p$ gives an
isomorphism between the ring $\Gamma$ defined in Section \ref{spsC}
and the ring $\Z[Q_1,Q_2,\ldots]$ of Schur $Q$-functions.
\end{example}

We wish to describe a tableau formula for $\Ti_\la(Z\,;X_k)$ which is
analogous to the expression (\ref{tabform}) for the Schur
functions. The following key observation is used to define the
relevant tableaux. In the Pieri rule (\ref{pruleC}) for the product
$u_p\cdot \Ti_\la(u)$, all partitions $\mu$ which appear on the right
hand side may be written as $\mu=(p+r,\nu)$ for some integer $r\geq 0$
and $k$-strict partition $\nu$ with $\nu\subset\la$. Moreover, if $p$
is sufficiently large (for example $p>|\la|+2k$) and we write
\begin{equation}
\label{pprod}
u_p\cdot \Ti_\la(u) = \sum_{r,\nu}2^{n(\la/\nu)}\,\Ti_{(p+r,\nu)}(u),
\end{equation}
with the sum over integers $r\geq 0$ and $k$-strict partitions
$\nu\subset\la$ with $|\nu| = |\la|-r$, then the $\nu$ which appear in
(\ref{pprod}) and the exponents $n(\la/\nu)$ are independent of
$p$. If this is the case, so that $\la\xrightarrow{p} (p+r,\nu)$ for
any $p>|\la|+2k$, then we say that $\la/\nu$ is a {\em $k$-horizontal
  strip}.

Let $\la$ and $\mu$ be any two $k$-strict partitions with $\mu\subset\la$.
A {\em $k$-tableau} $T$ of shape $\la/\mu$ is a sequence of $k$-strict
partitions
\[
\mu = \la^0\subset\la^1\subset\cdots\subset\la^r=\la
\]
such that $\la^i/\la^{i-1}$ is a $k$-horizontal strip for $1\leq i\leq
r$.  We represent $T$ by a filling of the boxes in $\la/\mu$ with
positive integers which is weakly increasing along each row and down
each column, such that for each $i$, the boxes in $T$ with entry $i$
form the skew diagram $\la^i/\la^{i-1}$. For any
$k$-tableau $T$ we define $n(T):=\sum_i n(\la^i/\la^{i-1})$ and let $c(T)$
denote the content vector of $T$.

Let {\bf P} denote the ordered alphabet
$\{1'<2'<\cdots<k'<1<2<\cdots\}$.  We say that the symbols
$1',\ldots,k'$ are {\em marked}, while the rest are {\em unmarked}.  A
{\em $k$-bitableau} $U$ of shape $\la$ is a filling of the boxes in
the diagram of $\la$ with elements of {\bf P} which is weakly
increasing along each row and down each column, such that the marked
entries are strictly increasing along each row, and the unmarked
entries form a $k$-tableau $T$. We define
\[
n(U):=n(T) \quad \text{and} \quad
(zx)^{c(U)}:= z^{c(T)} \,\prod_{j=1}^k x_j^{m_j} 
\]
where $m_j$ denotes the number of times that $j'$ appears in $U$.  For
any $k$-strict partition $\la$, we then have the {\em tableau formula}
\begin{equation}
\label{tableauxeq}
\Ti_\la(Z\,;X_k) = \sum_U 2^{n(U)}(zx)^{c(U)} 
\end{equation}
where the sum is over all $k$-bitableaux $U$ of shape $\la$. Using the
tableau formula (\ref{tabform}) for Schur polynomials, we can rewrite
equation (\ref{tableauxeq}) as 
\[
\Ti_\la(Z\,;X_k) = \sum_{\mu\subset\la}\sum_T 2^{n(T)}z^{c(T)} s_{\wt{\mu}}(X_k)
\]
with the sums over all partitions $\mu\subset\la$ and $k$-tableaux $T$
of shape $\la/\mu$, respectively.

\begin{example}
\label{Qex}
If $k=0$ and $\la$ is a strict partition, then (\ref{tableauxeq}) 
becomes
\begin{equation}
\label{Qeq}
Q_\la(Z) = \sum_T 2^{n(T)}z^{c(T)} 
\end{equation}
summed over all $0$-tableaux $T$ of shape $\la$. Equation (\ref{Qeq})
is a tableau formula for the Schur $Q$-functions.
\end{example}

\begin{example}
(a) For any integer $p\geq 0$, we have 
\[
\Ti_p(Z\,;X_k) = \sum_{j=0}^p \sum_{|\al|=p-j} 2^{\#\al}z^\al s_{1^j}(X_k)
=\sum_{j=0}^p Q_{p-j}(Z)e_j(X_k)
\]
where the second sum is over all compositions $\al$ with $|\al|=p-j$, 
and $\#\al$ denotes the number of indices $i$ such that $\al_i\neq 0$.

\medskip
\noin
(b) Assume that $k\geq 1$. Then we have
 \[
\Ti_{1^p}(Z\,;X_k) = \sum_{j=0}^p\sum_{|\al|=p-j} 2^{\#\al}z^\al s_j(X_k)
=\sum_{j=0}^p Q_{p-j}(Z)h_j(X_k).
\]
\end{example}

\medskip

The group $W_\infty$ is generated by the reflections $s_i$ for 
$i\geq 0$, and these generators are used to define reduced words and the 
length of signed permutations as in Section \ref{ac}.
The {\em   nilCoxeter algebra} $\cW_n$ of the hyperoctahedral group $W_n$ 
is the free associative algebra with unit
generated by the elements $\xi_0,\ldots,\xi_{n-1}$, modulo the relations
\[
\begin{array}{rclr}
\xi_i^2 & = & 0 & i\geq 0\ ; \\
\xi_i\xi_j & = & \xi_j\xi_i & |i-j|\geq 2\ ; \\
\xi_i\xi_{i+1}\xi_i & = & \xi_{i+1}\xi_i\xi_{i+1} & i\geq 1\ ; \\
\xi_0\xi_1\xi_0\xi_1 & = & \xi_1\xi_0\xi_1\xi_0.
\end{array}
\]
For any $w\in W_n$, choose a reduced word $a_1\cdots a_\ell$ for
$w$ and set $\xi_w := \xi_{a_1}\ldots \xi_{a_\ell}$. Then the $\xi_w$
are well defined and form a free $\Z$-basis of
$\cW_n$. We denote the coefficient of $\xi_w\in \cW_n$ in the
expansion of the element $\zeta\in \cW_n$ by $\langle \zeta,w\rangle$.

Let $t$ be an independent variable, define
\[
C(t) := (1+t \xi_{n-1})\cdots (1+t \xi_1)(1+t\xi_0)(1+t\xi_0)(1+t\xi_1) 
\cdots (1+t\xi_{n-1})
\]
and let $C(Z):=C(z_1)C(z_2)\cdots$. Choose an integer $k$ with $0\leq
k<n$ and set $A(X_k):=A(x_1)\cdots A(x_k)$. For any $w\in W_n$, the
(restricted) type C {\em mixed Stanley function} $J_w(Z\,;X_k)$ is
defined by
\[
J_w(Z\,;X_k) := \langle C(Z)A(X_k), w\rangle. 
\]
Clearly, $J_w$ is a power series in the $Z$ variables and a polynomial
in the $X_k$ variables, and has nonnegative integer coefficients. One
can show that $J_w$ is symmetric in the $Z$ and $X_k$ variables
separately.

When $w$ is the $k$-Grassmannian permutation associated to a
$k$-strict partition $\la$, then 
\[
J_w(Z\,;X_k)=\Ti_\la(Z\,;X_k).
\]
This equality may be generalized as follows.  We say that a signed
permutation $w=(w_1,\ldots,w_n)$ has a {\em descent} at position
$i\geq 1$ if $w_i>w_{i+1}$, and a descent at $i=0$ if and only if
$w_1<0$. A signed permutation $w\in W_n$ is {\em increasing up to $k$}
if it has no descents less than $k$. This condition is automatically
satisfied if $k=0$, and for positive $k$ it means that $0 < w_1 < w_2
< \cdots < w_k$. For any element $w\in W_n$ which is increasing up to
$k$, we have
\begin{equation}
\label{Jeq}
J_w(Z\,;X_k) = \sum_{\la}e^w_\la \,\Ti_\la(Z\,;X_k)
\end{equation}
summed over $k$-strict partitions $\la$ with $|\la|=\ell(w)$. The
integers $e^w_\la$ are nonnegative, that is, the function
$J_w(Z\,;X_k)$ is {\em theta positive} when $w$ is increasing up to
$k$. These coefficients have only one known combinatorial
interpretation, which is given below.

For positive integers $i<j$, define the reflections $t_{ij}$
and $\ov{t}_{ij},\ov{t}_{ii}$ in $W_\infty$ by their right actions
\begin{align*}
(w_1,\ldots,w_i,\ldots,w_j,\ldots)\,t_{ij} &= 
(w_1,\ldots,w_j,\ldots,w_i,\ldots), \\
(w_1,\ldots,w_i,\ldots,w_j,\ldots)\,\ov{t}_{ij} &= 
(w_1,\ldots,\ov{w}_j,\ldots,\ov{w}_i,\ldots), \ \ \mathrm{and} \\
(w_1,\ldots,w_i,\ldots)\,\ov{t}_{ii} &= 
(w_1,\ldots,\ov{w}_i,\ldots),
\end{align*} 
and let $\ov{t}_{ji} := \ov{t}_{ij}$.  For any $w\in W_\infty$ which
is increasing up to $k$, we construct a rooted tree $T^k(w)$ with root
$w$ and whose nodes are elements of $W_\infty$ as follows. Let $r$ be
the largest descent of $w$. If $w=1$ or $r=k$, then set
$T^k(w):=\{w\}$. Otherwise, let $s := \max(i>r\ |\ w_i < w_r)$ and
$\Phi(w) := \Phi_1(w)\cup \Phi_2(w)$, where
\begin{gather*}
\Phi_1(w) := \{wt_{rs}t_{ir}\ |\ 1\leq i < r \ \ \mathrm{and} \ \ 
\ell(wt_{rs}t_{ir}) = \ell(w) \}, \\
\Phi_2(w) := 
\{wt_{rs}\ov{t}_{ir}\ |\ i\geq 1 \ \ \mathrm{and} \ \
\ell(wt_{rs}\ov{t}_{ir}) = \ell(w) \}.
\end{gather*}
We define $T^k(w)$ recursively, by joining $w$ by an edge to each
$v\in \Phi(w)$, and attaching to each $v\in \Phi(w)$ its tree
$T^k(v)$. The finite tree $T^k(w)$ is the {\em $k$-transition tree} of
$w$, and its leaves are all $k$-Grassmannian elements of $W_\infty$.
The coefficient $e^w_\la$ in (\ref{Jeq}) is equal to the number of
leaves of $T^k(w)$ which have shape $\la$.

\begin{example}
Suppose that $k=1$. The $1$-transition tree for the signed permutation
$w=(3,\ov{1},2,5,4)$ in $W_5$ is shown below.
\medskip
\[
\includegraphics[scale=0.31]{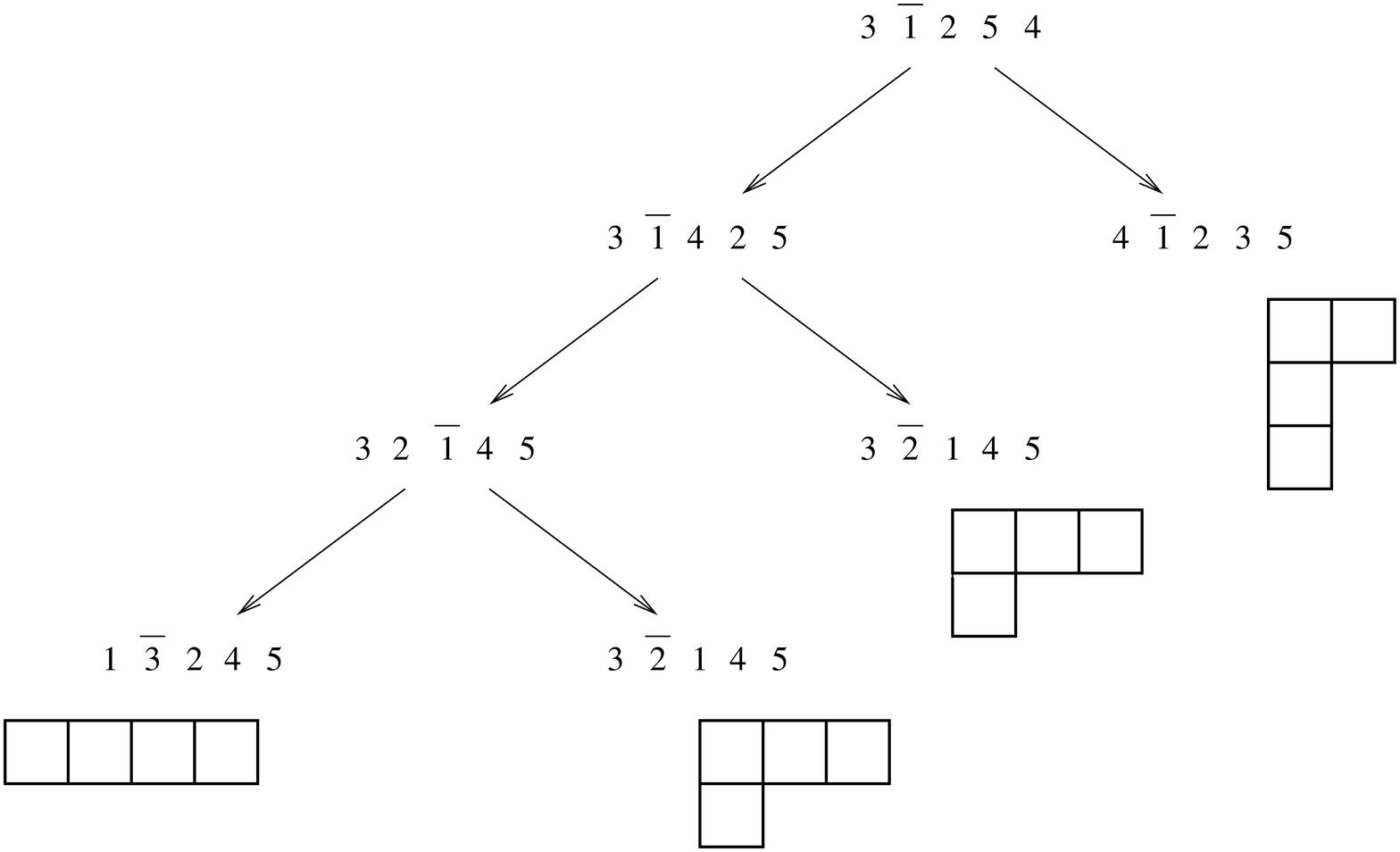}
\]
\medskip
\medskip
It follows that
\[
J_{3\ov{1}254}(Z\,;X_1)=\Ti_4(Z\,;X_1)+2\,\Ti_{3,1}(Z\,;X_1)+\Ti_{2,1,1}(Z\,;X_1).
\]
\end{example}

\subsection{Cohomology of flag manifolds}

Let $\{e_1,\ldots,e_{2n}\}$ denote the standard symplectic basis of
$E:=\C^{2n}$ and let $F_i=\langle e_1,\ldots, e_i\rangle$ be the
subspace spanned by the first $i$ vectors of this basis, as in Section
\ref{cGC}.  The group $G=\Sp_{2n}$ acts transitively on the space
of all complete isotropic flags in $E$, and the stabilizer of the flag
$F_\bull$ is a Borel subgroup $B$ of $G$. Let $T\subset B$ denote the
maximal torus of diagonal matrices in $G$, and the Weyl group
$W=N_G(T)/T\cong W_n$.

The parabolic subgroups $P$ of $\Sp_{2n}$ with $P\supset B$
correspond to sequences $a_1<\cdots < a_p$ of nonnegative 
integers with $a_p<n$. For any such $P$, the manifold $\X:=\Sp_{2n}/P$
parametrizes partial flags of subspaces
\[
E_\bull \ :\ 0 \subset E_p \subset \cdots \subset E_1 \subset E=\C^{2n}
\]
with $E_1$ isotropic and $\dim(E_r) = n-a_r$ for each $r\in [1,p]$. 
As usual, $E_r$ and $E$ will also denote the corresponding
tautological vector bundles over $\X$. The associated
parabolic subgroup $W_P$ of $W_n$ is generated by the
simple reflections $s_i$ for $i\notin\{a_1,\ldots, a_p\}$.

There is a canonical presentation of the cohomology ring of $\Sp_{2n}/B$,
which gives geometric significance to the variables which appear in
Section \ref{spsC}. Let $\IGam^{(n)}$ denote the ideal of $\Gamma[X_n]$
generated by the homogeneous elements of positive degree in
$\Gamma^{(n)}$, so that $\IGam^{(n)}=\langle
{}^nc_1,{}^nc_2,\ldots\rangle$. We then have ring isomorphism
\begin{equation}
\label{BorelC}
\HH^*(\Sp_{2n}/B,\Z) \cong \Gamma[X_n]/\IGam^{(n)}
\end{equation}
which maps each variable $c_p$ to $c_p(E/E_n)$ and $x_i$ to
$c_1(E_{n+1-i}/E_{n-i})$ for $1\leq i \leq n$. Furthermore, for any
parabolic subgroup $P$ of $\Sp_{2n}$, the projection map $\Sp_{2n}/B\to
\Sp_{2n}/P$ induces an injection $\HH^*(\Sp_{2n}/P,\Z)\hookrightarrow
\HH^*(\Sp_{2n}/B,\Z)$ on cohomology rings, and we have
\begin{equation}
\label{BorelPC}
\HH^*(\X,\Z) \cong \Gamma[X_n]^{W_P}/\IGam^{(n)}_P,
\end{equation}
where $\Gamma[X_n]^{W_P}$ denotes the $W_P$-invariant subring of
$\Gamma[X_n]$, and $\IGam_P^{(n)}$ is the ideal of $\Gamma[X_n]^{W_P}$ 
generated by ${}^nc_1,{}^nc_2,\ldots$

We have a decomposition
\[
\Sp_{2n} = \bigcup_{w \in W^P}Bw P
\]
where 
\[
W^P := \{w\in W_n\ |\ \ell(ws_i) = \ell(w)+1,\  \forall\, i \notin
\{a_1,\ldots, a_p\}, \ i<n\}
\]
is the set of minimal length $W_P$-coset representatives in $W_n$.
For each $w \in W^P$, the $B$-orbit of $wP$ in $\X$ is the {\em
  Schubert cell} $\Y^\circ_w:=Bw P/P$. The {\em Schubert variety}
$\Y_w$ is the closure of $\Y^\circ_w$ in $\Sp_{2n}/P$. Then
$\X_w:=\Y_{w_0w}$ has codimension $\ell(w)$ in $\X$, and its
cohomology class $[\X_w]$ is a {\em Schubert class}. The cell
decomposition of $\X$ implies that there is an isomorphism of abelian
groups
\[
\HH^*(\X,\Z) \cong \bigoplus_{w\in W^P}\Z[\X_w]
\]
which generalizes (\ref{IGrasseq}).

Recall that $E_r$ for $r\in [1,p]$ and $E$ denote the tautological and
trivial vector bundles over $\X$, of rank $n-a_r$ and $2n$,
respectively. For any $w\in W^P$, we then have
\begin{equation}
\label{Cpartfl}
[\X_w] = \sum_{\underline{\la}} f^w_{\underline{\la}}\,
\Ti^{(a_1)}_{\la^1}(E-E_1) s_{\la^2}(E_1-E_2)\cdots s_{\la^p}(E_{p-1}-E_p)
\end{equation}
in $\HH^*(\X,\Z)$, where the sum is over all sequences of
partitions $\underline{\la}=(\la^1,\ldots,\la^p)$ with $\la^1$ being 
$a_1$-strict, and the 
coefficients $f^w_{\underline{\la}}$ are given by 
\begin{equation}
\label{Ccoeff}
f^w_{\underline{\la}} := \sum_{u_1\cdots u_p = w}
e_{\la^1}^{u_1}c_{\la^2}^{u_2}\cdots c_{\la^p}^{u_p}
\end{equation}
summed over all factorizations $u_1\cdots u_p= w$ such that
$\ell(u_1)+\cdots + \ell(u_p)=\ell(w)$, $u_j\in S_\infty$ for $j\geq
2$, and $u_j(i)=i$ for all $j>1$ and $i \leq a_{j-1}$.  The
nonnegative integers $e^{u_1}_{\la^1}$ and $c^{u_i}_{\la^i}$ which
appear in the summands in (\ref{Ccoeff}) are the same as the ones in
equations (\ref{Jeq}) and (\ref{Geq}), respectively.  When $p=1$, the
partial flag manifold $\X$ is the isotropic Grassmannian
$\IG(n-a_1,2n)$, and formula (\ref{Cpartfl}) specializes to equation
(\ref{GiambC}). In general, the polynomial $\Ti^{(a_1)}_{\la^1}(E-E_1)$ 
in (\ref{Cpartfl}) is defined by pulling back the polynomial 
$\Ti^{(a_1)}_{\la^1}(c(E/E_1))$ under the natural projection map 
$\X \to \IG(n-a_1,2n)$ which sends a partial flag $E_\bull$ to $E_1$.

\begin{example}
\label{CSex}
Let $P=B$ be the Borel subgroup, so that the flag manifold
$\Sp_{2n}/B$ parametrizes flags of subspaces $0\subset E_n \subset
\cdots \subset E_1 \subset E=\C^{2n}$ with $E_1$ Lagrangian and
$\dim(E_i) = n+1-i$ for each $i\in [1,n]$. For each $i$, let
$x_i:=-c_1(E_i/E_{i+1})$, and observe that since $E/E_1\cong E_1^*$,
for any integer $p$, we have $c_p(E/E_1) = e_p(x_1,\ldots,x_n) =
e_p(X_n)$, using the definition of Chern classes. For every strict
partition $\la$, define the {\em $\wt{Q}$-polynomial}
$\wt{Q}_\la(X_n)$ by the formula
\[
\wt{Q}_\la(X_n):=\prod_{i<j}\frac{1-R_{ij}}{1+R_{ij}}\, e_\la(X_n).
\]
Since $a_1=0$, for any strict partition $\la$, we have
\[
\Ti^{(a_1)}_\la(E-E_1) = \prod_{i<j}\frac{1-R_{ij}}{1+R_{ij}}\,
c_\la(E/E_1) = \wt{Q}_\la(X_n).
\]
For any $w\in W_n$, we define the {\em symplectic Schubert
polynomial} $\CS_w(X_n)$ by 
\begin{equation}
\label{Teq}
\CS_w(X_n):= \sum_{v,\om,\la}e^v_\la\,\wt{Q}_\la(X_n)\AS_\om(-X_n)
\end{equation}
where the sum is over all factorizations $v\om = w$ and strict
partitions $\la$ with $\ell(v)+\ell(\om)=\ell(w)$, $\om\in S_n$, and
$|\la|=\ell(v)$. Employing computations similar to those in Example
\ref{ASex}, we then see that formula (\ref{Cpartfl}) is equivalent to
the statement that for any element $w\in W_n$, we have
$[\X_w]=\CS_w(X_n)$ in $\HH^*(\Sp_{2n}/B,\Z)$.
\end{example}

\section{Eta polynomials} 
\label{etapolys}

\subsection{Definition and Pieri rule}
\label{firstD}

In the theory of eta polynomials, we must distinguish between the  
case of level zero and that of positive level. Given any strict partition
$\la$, the {\em eta polynomial} $\Eta^{(0)}_\la(u)$ of level 
$0$ is defined by 
\[
\Eta^{(0)}_{\la}(u) := 2^{-\ell(\la)} \prod_{i<j}\frac{1-R_{ij}}{1+R_{ij}}\, u_{\la}.
\]
As a polynomial in the variables $u_p$, $\Eta^{(0)}_{\la}(u)$ may have
nonintegral coefficients. However, if we introduce new variables
$\bb_p$ such that $u_p=2\bb_p$ for each $p\geq 1$, then
$\Eta^{(0)}_{\la}$ is a polynomial in the $\bb_p$ with integer
coefficients.

Let $\B^{(0)}$ be the quotient of the polynomial ring 
$\Z[\bb_1,\bb_2,\ldots]$ modulo the ideal of relations
\[
\bb^2_p + 2\sum_{i=1}^{p-1}(-1)^i \bb_{p+i}\bb_{p-i}+(-1)^p\bb_{2p}
= 0
\ \ \ \text{for} \  p \geq 1.
\]
Then, for any $p\geq 0$, the Pieri rule
\begin{equation}
\label{thPieriD1}
\bb_p\cdot \Eta^{(0)}_\la(u) = \sum_{\mu} 2^{N'(\la,\mu)}\Eta^{(0)}_\mu(u)
\end{equation}
holds in $\B^{(0)}$, where the sum is over all strict partitions $\mu\supset
\la$ such that $|\mu|=|\la|+p$ and $\mu/\la$ is a horizontal strip,
and $N'(\la,\mu)$ is one less than the number of connected components of
$\mu/\la$.

Assume next that $k \geq 1$, and let
$\bb_1,\ldots,\bb_{k-1},\bb_k,\bb_k',\bb_{k+1},\ldots$ be independent
variables related to the variables $u_1,u_2,\ldots$ by the equations
\begin{equation}
\label{utob}
u_p=
\begin{cases}
\bb_p &\text{if $p< k$},\\
\bb_k+\bb_k' &\text{if $p=k$},\\
2\bb_p &\text{if $p> k$}.
\end{cases}
\end{equation}
The eta polynomials $\Eta^{(k)}_{\la}(u)$ lie in the ring
$\Z[\bb_1,\ldots,\bb_{k-1},\bb_k,\bb_k',\bb_{k+1},\ldots]$, and are indexed
by {\em typed} $k$-strict partitions $\la$.

A {\em typed $k$-strict partition} is a pair consisting of a $k$-strict
partition $\la$ together with an integer in $\{0,1,2\}$ called the
{\em type} of $\la$, and denoted $\type(\la)$, such that
$\type(\la)=0$ if and only if $\la_i\neq k$ for all $i\geq 1$. The type
is usually omitted from the notation for the pair $(\la,\type(\la))$.

For a general typed $k$-strict
partition $\la$, we define the operator
\begin{equation}
\label{raiseopD}
R^{\la} := \prod_{i<j} (1-R_{ij})\prod_{\la_i+\la_j \geq 2k+j-i}
(1+R_{ij})^{-1}
\end{equation}
where the first product is over all pairs $i<j$ and second product is
over pairs $i<j$ such that $\la_i+\la_j \geq 2k+j-i$.  Let $R$ be any
finite monomial in the operators $R_{ij}$ which appears in the
expansion of the formal power series $R^\la$ in (\ref{raiseopD}).  If
$\type(\la)=0$, then set $R \star u_{\la} := u_{R \,\la}$. Suppose
that $\type(\la)\neq 0$, let $r$ be the least index such that
$\la_r=k$, and set $$\wh{\al} :=
(\al_1,\ldots,\al_{r-1},\al_{r+1},\ldots,\al_\ell)$$ for any integer
sequence $\al$ of length $\ell$. If $R$ involves any factors $R_{ij}$
with $i=r$ or $j=r$, then let $R \star u_{\la} := \frac{1}{2}\,u_{R
  \,\la}$. If $R$ has no such factors, then let
\[
R \star u_{\la} := \begin{cases}
\bb_k \,u_{\wh{R \,\la}} & \text{if  $\,\type(\la) = 1$}, \\
\bb'_k \, u_{\wh{R \,\la}} & \text{if  $\,\type(\la) = 2$}.
\end{cases}
\]

We define the {\em eta polynomial} $\Eta^{(k)}_\la(u)$ of level $k$ by
\[
\Eta^{(k)}_\la(u) := 2^{-\ell_k(\la)} R^{\la} \star u_\la.
\]
Here the {\em $k$-length} $\ell_k(\la)$ of a (typed) $k$-strict
partition $\la$ is the number of parts $\la_i$ which are strictly
greater than $k$.  It is easy see that $\Eta_\la^{(k)}(u)$ is a
polynomial in the variables $\bb_p$ and $\bb'_k$ with integer
coefficients. We will write $\Eta_\la(u)$ for $\Eta^{(k)}_\la(u)$ when
the level $k$ is understood.

\begin{example}
(a) Consider the typed $2$-strict partition 
$\la=(3,2,2)$ with $\type(\la)=2$. Then we have
\begin{align*}
\Eta^{(2)}_\la(u) & = \frac{1}{2}
\frac{1-R_{12}}{1+R_{12}}(1-R_{13})(1-R_{23}) \star u_{3,2,2} \\
&= \frac{1}{2}(1-2R_{12}+2R_{12}^2 - 2 R_{12}^3)(1-R_{13}-R_{23}+R_{13}R_{23}) 
\star u_{3,2,2} \\
&= \bb_3\bb'_2(\bb_2+\bb_2') - \bb_3^2\bb_1 + \bb_4\bb_3 - \bb_4\bb'_2\bb_1 
+\bb_6\bb_1-\bb_7.
\end{align*}

\noin
(b) Let $\la$ be a $k$-strict partition with $\la_i=k\geq 1$ for some $i$, 
and let $\Eta_\la(u)$ and $\Eta'_\la(u)$ denote the eta polynomials of level
$k$ indexed by $\la$ of type 1 and 2, respectively. Then we have
\[
\Eta_\la(u)+ \Eta'_\la(u) = 2^{-\ell_k(\la)} R^{\la} \, u_\la
\]
where $R^\la$ denotes the operator (\ref{raiseopD}).
\end{example}

Let $\B^{(k)}$ be the quotient
of polynomial ring  $\Z[\bb_1,\ldots,\bb_{k-1},\bb_k,\bb_k',\bb_{k+1},\ldots]$
modulo the ideal of relations
\begin{gather}
\label{stabrelD1}
\bb_p^2+ \sum_{i=1}^p(-1)^i\bb_{p+i}u_{p-i} = 0 \ \ \ \text{for
$p>k$}, \\
\label{stabrelD2}
\bb_k\bb'_k+\sum_{i=1}^k (-1)^i\bb_{k+i}\bb_{k-i} = 0,
\end{gather}
where the $u_i$ obey the equations (\ref{utob}). For every typed
$k$-strict partition $\la$, we define a monomial $\bb_\la$ as
follows. If $\type(\la)\neq 2$, then set
$\bb_\la:=\bb_{\la_1}\bb_{\la_2}\cdots$. If $\type(\la)= 2$ then
define $\bb_\la$ by the same product formula, but replacing each
occurrence of $\bb_k$ with $\bb'_k$. The monomials $\bb_\la$ and the
polynomials $\Eta_\la(u)$ as $\la$ runs over all typed
$k$-strict partitions form two $\Z$-bases of the graded ring
$\B^{(k)}$. The Pieri rule for the products $\bb_p\cdot\Eta_\la(u)$ 
holds only modulo the relations (\ref{stabrelD1}) and (\ref{stabrelD2});
again we need some further definitions to state it.

We say that the box in row $r$ and column $c$ of a $k$-strict
partition $\lambda$ is {\em $k'$-related} to the box in row $r'$ and
column $c'$ if $|c-k-1/2|+r = |c'-k-1/2|+r'$.  For example, the two
grey boxes in the following partition are $k'$-related.
\[ \pic{0.65}{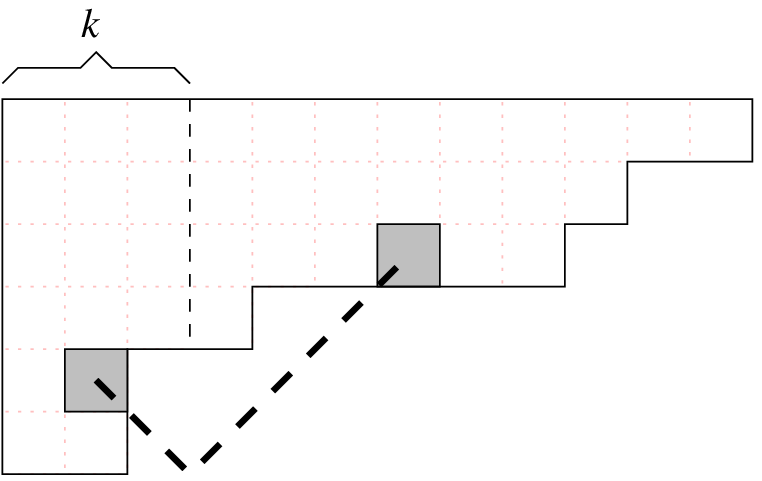} \]
For any two $k$-strict partitions $\lambda$ and $\mu$, the
relation $\lambda \xrightarrow{p} \mu$ is defined as in Section
\ref{firstC}, but replacing `$k$-related' by `$k'$-related'
throughout. The set $\D$ of boxes of $\mu\ssm \la$
is defined in the same way, and the integer $N'(\la,\mu)$ is equal to
the number (respectively, one less than the number) of connected
components of $\D$, if $p\leq k$ (respectively, if $p>k$).

If $\lambda$ and $\mu$ are typed $k$-strict partitions, then we write
$\lambda \xrightarrow{p} \mu$ if the underlying $k$-strict partitions
satisfy $\lambda \xrightarrow{p} \mu$, with the added condition that
$\type(\la)+\type(\mu)\neq 3$.  Let $c(\la,\mu)$ be the number of
columns of $\mu$ among the first $k$ which do not have more boxes than
the corresponding column of $\la$, and
\[
d(\la,\mu) := c(\la,\mu)+\max(\type(\la),\type(\mu)).
\]
If $p\neq k$, then set $\delta_{\la\mu}=1$. If $p=k$ and 
$N'(\la, \mu)>0$, then set
\[
\delta_{\la\mu}=\delta'_{\la \mu}:=1/2,
\]
while if $N'(\la, \mu)=0$, define
\[
\delta_{\la \mu} := \begin{cases} 
1 & \text{if $d(\la,\mu)$ is odd}, \\
0 & \text{otherwise}
\end{cases}
\qquad \mathrm{and}
\qquad 
\delta'_{\la\mu} = \begin{cases}
1 & \text{if $d(\la,\mu)$ is even}, \\
0 & \text{otherwise.}
\end{cases}
\]

For any typed $k$-strict partition $\la$ and $p\geq 0$, we then
have the {\em Pieri rule}
\begin{equation}
\label{pruleD}
\bb_p \cdot\Eta_{\lambda}(u) = \sum_\mu
\delta_{\la\mu}\,2^{N'(\la,\mu)}\,\Eta_\mu(u),
\end{equation}
in $\B^{(k)}$, where the sum over all typed $k$-strict partitions $\mu$ 
such that $\la\xrightarrow{p}\mu$. Furthermore, the product
$\bb'_k\cdot\Eta_{\lambda}(u)$ is obtained by replacing $\delta_{\lambda
\mu}$ with $\delta'_{\lambda \mu}$ throughout.

\subsection{Cohomology of Grassmannians} 
\label{cGD}

Equip the vector space $\C^{2n}$ with the nondegenerate symmetric
bilinear form $(\ ,\,)$ defined by the conditions $(e_i,e_j)=0$ for
$i+j\neq 2n+1$ and $(e_i,e_{2n+1-i})=1$ for $1\leq i \leq n$. The
special orthogonal group $\SO_{2n}(\C)$ is the subgroup of
$\SL_{2n}(\C)$ consisting of those elements $g$ such that
$(gv_1,gv_2)=(v_1,v_2)$, for every $v_1,v_2\in \C^{2n}$. We say that a
subspace $V$ is {\em isotropic} if the restriction of $(\ ,\,)$
to $V$ vanishes identically. Since the form is nondegenerate, we
have $\dim(V)\leq n$ for any isotropic subspace $V$. For
each $i\in [1,2n]$, let $F_i$ denote the $\C$-linear span of
$e_1,\ldots,e_i$.

Fix an integer $k$ with $0\leq k \leq n-1$. If $k\geq 1$, then the
orthogonal Grassmannian $\OG=\OG(n-k,2n)$ parametrizes isotropic
linear subspaces of dimension $n-k$ in $\C^{2n}$. The group $\SO_{2n}$
acts transitively on $\OG(n-k,2n)$, and we have $\OG(n-k,2n) =
\SO_{2n}/P_k$, where $P_k$ is a maximal parabolic subgroup of
$\SO_{2n}$. When $k=0$, the locus of maximal isotropic subspaces has
two isomorphic connected components, called the two {\em families},
each of which is a single $\SO_{2n}$-orbit.  The orthogonal
Grassmannian $\OG(n,2n)=\SO_{2n}/P_0$ parametrizes one of these
components, which we take to be the family containing $F_n$.

We agree that when $k=0$, a typed $0$-strict partition is the same as
a strict partition, and that all such partitions have type 1. The
Schubert cells $\X^\circ_\la$ in $\OG(n-k,2n)$ are indexed by the
typed $k$-strict partitions $\la$ whose diagrams are contained in an
$(n-k)\times (n+k-1)$ rectangle. We have
\[
X^\circ_\la := \{ V \in \OG \mid \dim(V \cap
   F_r) = \#\{j\ |\ p_j(\la)\leq r\} \ \ \ \forall\, r \},
\]
where the strictly increasing index function
$\{p_j(\lambda)\}_{1\leq j \leq n-k}$ is defined by
\begin{multline*} p_j(\lambda) := n+k+j-\lambda_j - 
   \#\{\,i<j\ |\ \la_i+\la_j \geq 2k+j-i\,\} \\
   {} - \begin{cases} 
      1 & \text{if $\lambda_j > k$, or $\lambda_j=k < \lambda_{j-1}$ and
        $n+j+\type(\lambda)$ is odd}, \\
      0 & \text{otherwise}.
   \end{cases}
\end{multline*}
The Schubert variety $\X_\la$ is the closure of the
$\X^\circ_\la$, and has codimension $|\la|$ in $\OG$. 
If $[\X_\la]$ denotes the cohomology class of $\X_\la$ in
$\HH^{2|\la|}(\OG,\Z)$, then we have a group isomorphism
\begin{equation}
\label{OGrasseq}
\HH^*(\OG(n-k,2n),\Z) \cong \bigoplus_\la \Z [\X_\la].
\end{equation}

The varieties $\X_p$ for $1\leq p \leq n+k-1$, together with $\X_k'$
when $k\geq 1$, are the {\em special Schubert varieties}, and their
classes in cohomology are the {\em special Schubert classes}.  The
convention here is that $\X_k$ (respectively, $\X_k'$) is indexed by
the partition $k$ of type 1 (respectively, type 2). As in the Lie
types A and C, the special Schubert varieties $\X_p$ and $\X_k'$ can
be viewed as the locus of all isotropic linear subspaces $V$ which meet
a given isotropic or coisotropic linear subspace nontrivially.

Let $Q\to \OG(n-k,2n)$ denote the 
universal quotient vector bundle over $\OG$, of rank $n+k$. For $k=0$,
we have $c_p(Q)=2[\X_p]$ for all $p\geq 1$, while for $k\geq 1$, we have
\[
c_p(Q)= \begin{cases} [\X_p] & \text{if $p<k$}, \\
[\X_k]+[\X_k'] & \text{if $p=k$}, \\
2[\X_p] & \text{if $p>k$}
\end{cases}
\]
in $\HH^{2p}(\OG(n-k,2n),\Z)$, in agreement with (\ref{utob}). We can
now state the {\em Giambelli formula}
\begin{equation}
\label{GiambD}
[\X_\la] = \Eta_{\la}(c(Q))
\end{equation}
where the polynomial $\Eta_{\la}(c(Q))$ is obtained from
$\Eta^{(k)}_{\la}(u)$ by performing the substitutions $\bb_p\mapsto [\X_p]$
and $\bb'_k\mapsto [\X'_k]$ for every integer $p$.

Furthermore, the {\em Pieri rules} (\ref{thPieriD1}) and (\ref{pruleD}) hold
in $\HH^*(\OG(n-k,2n),\Z)$. For instance, the latter rule is valid for
$k\geq 1$ and states that
\begin{equation}
\label{pruleD2}
[\X_p]\cdot [\X_\la] = \sum_{\mu} \delta_{\la\mu}\,2^{N'(\la,\mu)}[\X_\mu]
\end{equation}
summed over all typed $k$-strict partitions 
$\mu$ such that $\la\xrightarrow{p}\mu$ and the diagram of $\mu$ fits 
in an $(n-k)\times (n+k-1)$ rectangle. Moreover, the product 
$[\X'_k]\cdot [\X_\la]$ is obtained by replacing $\delta_{\lambda
\mu}$ with $\delta'_{\lambda \mu}$ in (\ref{pruleD2}).

\begin{example}
For the Grassmannian $\OG(5,14)$ we have $n=7$ and $k=2$. Let $\la$
denote the partition $(8,7,2,1,1)$ of type 1. We then have the 
Pieri formulas
\begin{align*}
[\X_2]\cdot [\X_\la]  &= [\X_{8,7,4,1,1}]+[\X_{8,7,3,2,1}]+[\X_{8,7,6}] \\ 
[\X'_2]\cdot [\X_\la] &= [\X_{8,7,4,1,1}]+[\X_{8,7,3,2,1}]
\end{align*}
where the indexing partitions on the right hand side are all of type 0 or 1.
\end{example}

The ring $\B^{(k)}$ is naturally isomorphic to the stable cohomology
ring
\[
\IH(\OG_k) = \lim_{\longleftarrow}\HH^*(\OG(n-k,2n),\Z)
\]
of the orthogonal Grassmannian $\OG$, where the inverse limit is
defined as in Section \ref{cGC}. Under this isomorphism, the variables
$\bb_p$ and $\bb'_k$ map to the special Schubert classes $[\X_p]$ and
$[\X'_k]$ in the cohomology ring of $\OG$.

\subsection{Symmetric polynomials}
\label{spsD}

Let $b:=(b_1,b_2,\ldots)$ be a sequence of commuting variables, and
set $b_0:=1$ and $b_p=0$ for $p<0$. Consider the graded ring $\Gamma'$
which is the quotient of the ring $\Z[b]$ modulo the ideal
generated by the relations
\[
b_p^2+2\sum_{i=1}^{p-1}(-1)^i b_{p+i}b_{p-i}+(-1)^p b_{2p}=0, \ \ \ 
\text{for all $p\geq 1$}.
\]
The ring $\Gamma'$ is isomorphic to $\B^{(0)}$ and to the stable 
cohomology ring
\[
\lim_{\longleftarrow}\HH^*(\OG(n,2n),\Z)
\]
of the maximal orthogonal Grassmannian $\OG(n,2n)$, with the variables
$b_p$ mapping to the special Schubert classes $[\X_p]$. We regard $\Gamma$
as a subring of $\Gamma'$ using the injection which sends $c_p$ to 
$2b_p$ for all $p\geq 1$.

The Weyl group $\wt{W}_n$ for the root system $\text{D}_n$ is the
subgroup of $W_n$ consisting of all signed permutations with an even
number of sign changes.  The group $\wt{W}_n$ is an extension of $S_n$
by the element $s_\Box=s_0s_1s_0$, which acts on the right by
\[
(w_1,w_2,\ldots,w_n)s_\Box=(\ov{w}_2,\ov{w}_1,w_3,\ldots,w_n).
\]

Fix $n\geq 2$ and let $X_n:=(x_1,\ldots,x_n)$. There is a natural action of 
$\wt{W}_n$ on $\Gamma'[X_n]$ which extends the 
action of $S_n$ on $\Z[X_n]$, defined as follows. 
The simple reflections $s_i$ for $i>0$ act
by interchanging $x_i$ and $x_{i+1}$ and leaving all the remaining
variables fixed. The reflection $s_\Box$ maps 
$(x_1,x_2)$ to $(-x_2,-x_1)$, fixes the $x_j$ for $j\geq 3$, 
and satisfies, for any $p\geq 1$,
\begin{equation}
\label{sBx}
s_\Box(b_p) :=
b_p+(x_1+x_2)\sum_{j=0}^{p-1}\left(\sum_{a+b=j}x_1^ax_2^b\right)
c_{p-1-j}.
\end{equation}
If $t$ is a formal variable which is fixed by $s_\Box$, then we 
express equation (\ref{sBx}) using generating functions as
\[
s_\Box\left(\sum_{p=0}^\infty c_pt^p\right) = \frac{1+x_1t}{1-x_1t}\cdot
\frac{1+x_2t}{1-x_2t}\cdot\left(\sum_{p=0}^\infty c_pt^p\right).
\]

For every integer $p$, define an element ${}^nb_p$ of $\Gamma'[X_n]$ by
\[
{}^nb_p:=\begin{cases} 
e_p(X_n) + 2\sum_{i=0}^{p-1} e_i(X_n)b_{p-i} & \text{if $p<n$}, \\
\sum_{i=0}^p e_i(X_n)b_{p-i} & \text{if $p\geq n$},
\end{cases}
\]
let
\[
{}^nb'_n:=\sum_{i=0}^{n-1}e_i(X_n)b_{n-i}
\]
and set $B^{(n)}:= \Z[{}^nb_1, \ldots, {}^nb_{n-1},{}^nb_n, {}^nb'_n,
{}^nb_{n+1}, \ldots]$. Observe that we have
\[
{}^nc_p=\begin{cases}
{}^nb_p & \text{if $p<n$}, \\
{}^nb_n+{}^nb'_n & \text{if $p=n$}, \\
2\cdot {}^nb_p & \text{if $p>n$}
\end{cases}
\]
and thus $\Gamma^{(n)}$ is a subring of $B^{(n)}$. 

Let $\Gamma'[X_n]^{\wt{W}_n}$ denote the subring of
$\wt{W}_n$-invariants in $\Gamma'[X_n]$. Then there is an equality
\begin{equation}
\label{Binv}
\Gamma'[X_n]^{\wt{W}_n}=B^{(n)}= \Z[{}^nb_1, \ldots, {}^nb_{n-1},{}^nb_n, {}^nb'_n,
{}^nb_{n+1}, \ldots].
\end{equation}

The map which sends $\bb_p$ to ${}^nb_p$ for every integer $p$ and $\bb'_n$
to ${}^nb'_n$ induces a ring 
isomorphism $\B^{(n)}\cong B^{(n)}$. We therefore have
\begin{equation}
\label{GTieqD}
B^{(n)}=\bigoplus_\la \Z\, \Eta_\la(X_n)
\end{equation}
where the sum is over all typed $n$-strict partitions $\la$, and the
polynomial $\Eta_\la(X_n)$ is obtained from the eta polynomial
$\Eta^{(n)}_\la(u)$ by making the substitutions $\bb_p\mapsto {}^nb_p$ for
all $p$ and $\bb'_n\mapsto {}^nb'_n$. In other words, we have
\[
\Eta_\la(X_n)=\Eta^{(n)}_\la(X_n) := 2^{-\ell_n(\la)} R^\la \star ({}^nc)_\la.
\]

For each $r\geq 2$, we embed $\wt{W}_r$ in $\wt{W}_{r+1}$ by adding
the element $r+1$ which is fixed by $\wt{W}_r$, and set $\wt{W}_\infty
:=\cup_r \wt{W}_r$.  Let $w$ be a signed permutation in
$\wt{W}_\infty$.  Define a strict partition $\mu(w)$ whose parts are
the absolute values of the negative entries of $w$ minus one, arranged
in decreasing order.  Let the {\em A-code} of $w$ be the sequence
$\gamma$ with $\gamma_i:=\#\{j>i\ |\ w_j<w_i\}$, define a
partition $\delta(w)$ whose parts are the nonzero entries $\gamma_i$
arranged in weakly decreasing order, and let $\nu(w)$ be the conjugate
of $\delta(w)$. The {\em shape} of $w$ is defined to be the partition
$\la(w):=\mu(w)+\nu(w)$.

\begin{example}
(a) An element $w$ of $\wt{W}_\infty$ is $n$-Grassmannian if
  $\ell(ws_i)>\ell(w)$ for all $i\neq n$. The type of an
  $n$-Grassmannian element $w$ is 0 if $|w_1|=1$, and 1 (respectively,
  2) if $w_1>1$ (respectively, if $w_1<-1$). There is a type
  preserving bijection between the $n$-Grassmannian elements of
  $\wt{W}_\infty$ and typed $n$-strict partitions, given as follows.
  If the element $w$ corresponds to the typed $n$-strict partition
  $\la$, then for each $j\geq 1$, we have
\[
\la_i=\begin{cases} 
n-1+|w_{n+i}| & \text{if $w_{n+i}<0$}, \\
\#\{r\leq n\, :\, |w_r|> w_{n+i}\} & \text{if $w_{n+i}>0$}.
\end{cases}
\]
The shape $\la(w)$ of $w$ agrees with the typed $n$-strict partition 
associated to $w$, if $w$ has type 0 or 1. However, this may fail if
$w_1<-1$, for instance the $2$-Grassmannian element 
$v:=\ov{3}5\ov{1}24$ in $\wt{W}_5$ is associated to the
typed partition of shape $(2,2,1)$, while $\la(v)=(3,1,1)$.

\medskip
\noin
(b) The longest element of $\wt{W}_n$ is given by
\[
\wt{w}_0=\left\{ \begin{array}{cl}
           (\ov{1},\ldots,\ov{n}) & \mathrm{ if } \ n \ \mathrm{is} \ 
             \mathrm{even}, \\
           (1,\ov{2},\ldots,\ov{n}) & \mathrm{ if } \ n \ \mathrm{is} \ 
             \mathrm{odd}.
             \end{array} \right.
\]
Then $\mu(\wt{w}_0)=\delta_{n-1}$, $\nu(\wt{w}_0)=
\delta_{n-1}$, and $\la(\wt{w}_0)=2\delta_{n-1}=(2n-2,2n-4,\ldots,
2)$.
\end{example}

Let $w$ be an $n$-Grassmannian element of $\wt{W}_\infty$ with corresponding 
partition $\lambda(w)$, and let $\wt{w}_0$ be the longest element of $\wt{W}_n$.
Define
\[
{}^{\nu(w\wt{w}_0)}P_{\la(w\wt{w}_0)} := 2^{-r}\,\prod_{i<j} \frac{1-R_{ij}}{1+R_{ij}}
\, {}^{\nu(w\wt{w}_0)}c_{\la(w\wt{w}_0)},
\]
where $r$ is the length of the partition $\la(w\wt{w}_0)$.
The alternating operator $\cA''$ on $\Gamma'[X_n]$ is given by
\[
\cA''(f):=\sum_{w\in \wt{W}_n} (-1)^{\ell(w)}w(f).
\]
We then have
\begin{equation}
\label{mScheqD}
\Eta_{\la(w)}(X_n) 
=(-1)^{n(n-1)/2}\cdot 2^{n-1} \left. 
\cA''\left({}^{\nu(w\wt{w}_0)}P_{\la(w\wt{w}_0)}\right)\right\slash
\cA''\left(x^{\la(\wt{w}_0)}\right)
\end{equation}
in $\Gamma'[X_n]$.

\subsection{Algebraic combinatorics}
\label{acD}

Let $k$ be a nonnegative integer. In this section, we will define the 
formal power series $\Eta^{(k)}_\la(Z\,;X_k)$ which are the analogues in Lie
type D of the power series $\Ti^{(k)}_\la(Z\,;X_k)$ discussed in Section 
\ref{acC}. Their definition is easy when $k=0$: in this case, the index
$\la$ is a strict partition and we have
\[
\Eta^{(0)}_\la(Z) := 2^{-\ell(\la)} \Ti^{(0)}_\la(Z) = P_\la(Z),
\]
where $P_\la(Z)$ is a Schur $P$-function. The map which sends $\bb_j$
to $P_j:=P_j(Z)=Q_j(Z)/2$ for every positive integer $j$ gives an
isomorphism between the ring $\B^{(0)}\cong\Gamma'$ defined in Section
\ref{spsD} and the ring $\Z[P_1,P_2,\ldots]$ of Schur
$P$-functions. Since $\Eta^{(0)}_\la(Z)$ is a scalar multiple of
$\Ti^{(0)}_\la(Z)=Q_\la(Z)$, we immediately obtain a tableau formula
for $\Eta^{(0)}_\la(Z)$ from the tableau formula (\ref{Qeq}) for the
Schur $Q$-functions.

Assume next that $k\geq 1$. Set
\[
\eta_r=\eta_r(Z\,;X_k) := \begin{cases}
e_r(X_k) + 2\sum_{i=0}^{r-1} P_{r-i}(Z)e_i(X_k) & \text{if $r<k$}, \\
 \sum_{i=0}^r P_{r-i}(Z) e_i(X_k) & \text{if $r\geq k$}
\end{cases}
\]
and $$\eta'_k=\eta'_k(Z\,;X_k) =  \sum_{i=0}^{k-1} P_{k-i}(Z) e_i(X_k).$$ 
For any $r\geq 0$, if $\ti_r$ is defined by equation (\ref{tieq}),
then we have
\[
\ti_r=
\begin{cases}
\eta_r &\text{if $r< k$},\\
\eta_k+\eta_k' &\text{if $r=k$},\\
2\eta_r &\text{if $r> k$}.
\end{cases}
\]

By definition, for every typed $k$-strict partition $\la$,
the eta polynomial $\Eta_\la(Z\,;X_k)$ is obtained from
$\Eta^{(k)}_\la(u)$ by making the substitutions $\bb_r\mapsto \eta_r$
for every integer $r$ and $\bb'_k\mapsto \eta'_k$. In other words, we have
\[
\Eta_\la(Z\,;X_k):= 2^{-\ell_k(\la)} \,R^\la\star \ti_\la.
\]
We proceed to give a tableau formula for $\Eta_\la(Z\,;X_k)$ which is 
analogous to the formulas (\ref{tabform}) and (\ref{tableauxeq}).

Let $\la$ and $\mu$ be $k$-strict partitions with $\mu\subset\la$,
and choose any $p> |\la|+2k-1$. If
$|\la|=|\mu|+r$ and $\la\xrightarrow{p}(p+r,\mu)$, then we 
say that $\la/\mu$ is a {\em $k'$-horizontal strip}.  
We call a box in row $r$ and column $c$ of a Young diagram 
a {\em left box} if $c \leq k$ and a {\em
right box} if $c>k$.
If $\mu\subset\la$ are two $k$-strict
partitions such that $\la/\mu$ is a $k'$-horizontal strip,
we define $\la_0=\mu_0=+\infty$ and agree that the
diagrams of $\la$ and $\mu$ include all boxes $[0,c]$ in row zero.
We let $\E$ denote the set of right boxes of $\mu$
(including boxes in row zero) which are bottom boxes of $\la$ in their
column and are not $(k-1)$-related to a left box of
$\la/\mu$.

If $\la$ and $\mu$ are typed $k$-strict partitions with
$\mu\subset\la$, we say that $\la/\mu$ is a
{\em typed $k'$-horizontal strip} if the underlying $k$-strict
partitions are such that $\la/\mu$ is a $k'$-horizontal strip and in
addition $\type(\la)+\type(\mu)\neq 3$. In this case we let
$n(\la/\mu)$ denote  the number of connected components of
$\E$ minus one.

Suppose that $\la$ is any typed $k$-strict partition.  Let {\bf P$'$}
denote the ordered alphabet
$\{\wh{1}<\wh{2}<\cdots<\wh{k}<1,1^\circ<2,2^\circ<\cdots\}$.  We say 
that the symbols $\wh{1},\ldots,\wh{k}$ are {\em marked}, while the
rest are {\em unmarked}. A {\em typed $k'$-tableau} $T$ of shape
$\la/\mu$ is a sequence of typed $k$-strict partitions
\[
\mu = \la^0\subset\la^1\subset\cdots\subset\la^r=\la
\]
such that $\la^i/\la^{i-1}$ is a typed $k'$-horizontal strip for
$1\leq i\leq r$.  We represent $T$ by a filling of the boxes in
$\la/\mu$ with unmarked elements of {\bf P$'$} which is weakly
increasing along each row and down each column, such that for each
$i$, the boxes in $T$ with entry $i$ or $i^\circ$ form the skew
diagram $\la^i/\la^{i-1}$, and we use $i$ (resp.\ $i^\circ$) if and
only if $\type(\la^i)\neq 2$ (resp.\ $\type(\la^i)=2$), for every
$i\geq 1$.  For any typed $k'$-tableau $T$ we define $n(T)=\sum_i
n(\la^i/\la^{i-1})$ and let $c(T)=(r_1, r_2, \ldots)$ be the content
vector of $T$, so that $r_i$ denotes the number of times that $i$ or
$i^\circ$ appears in $T$, for each $i\geq 1$.

A {\em typed $k'$-bitableau} $U$ of shape $\la$ is a filling of the
boxes in the diagram of $\la$ with elements of {\bf P$'$} which is
weakly increasing along each row and down each column, such that the
unmarked entries form a typed $k'$-tableau $T$ of shape $\la/\mu$ with
$\type(\mu)\neq 2$, and the marked entries are a filling of $\mu$
which is strictly increasing along each row. We define
\[
n(U)=n(T) \quad \text{and} \quad
(zx)^{c(U)}= z^{c(T)} \,\prod_{j=1}^k x_j^{m_j} 
\]
where $m_j$ denotes the number of times that $\wh{j}$ appears in $U$.
For any typed $k$-strict partition $\la$, we then have the {\em
  tableau formula}
\begin{equation}
\label{tableauxeqD}
\Eta_\la(Z\,;X_k) = \sum_U 2^{n(U)}(zx)^{c(U)} 
\end{equation}
where the sum is over all typed $k'$-bitableaux $U$ of shape $\la$. 
Using the tableau formula (\ref{tabform}), we can rewrite
equation (\ref{tableauxeqD}) as 
\begin{equation}
\label{2ndeq}
\Eta_\la(Z\,;X_k) = \sum_{\mu\subset\la}\sum_T 2^{n(T)}z^{c(T)} s_{\wt{\mu}}(X_k)
\end{equation}
with the sums over all partitions $\mu\subset\la$ and typed $k'$-tableaux $T$
of shape $\la/\mu$, respectively.

\begin{example}
Suppose that $k=1$, and for any $r\geq 1$, let $\Eta_{1^r}(Z\,;X_1)$
and $\Eta'_{1^r}(Z\,;X_1)$ denote the eta polynomials indexed by $1^r$ of
type $1$ and $2$, respectively. Then for any integer $r\geq 0$, we
deduce from equation (\ref{2ndeq}) that
\begin{align*}
\Eta_r(Z\,;X_1) &= P_r(Z)+P_{r-1}(Z)x_1,\\
\Eta_{1^r}(Z\,;X_1) &= P_r(Z) + 2 P_{r-1}(Z)x_1 + \cdots + 
2 P_1(Z)x_1^{r-1}+x_1^r, \ \text{and}\\
\Eta'_{1^r}(Z\,;X_1) &= P_r(Z).
\end{align*}
\end{example}

\medskip

Consider the set $\N_\Box :=\{\Box,1,\ldots\}$ whose members index the
simple reflections in $\wt{W}_\infty$. These elements generate the
group $\wt{W}_\infty$ and are used to define reduced words and the
length of signed permutations as in Section \ref{acC}.  The {\em
  nilCoxeter algebra} $\wt{\cW}_n$ of the group $\wt{W}_n$ is the free
associative algebra with unit generated by the elements
$\xi_\Box,\xi_1,\ldots,\xi_{n-1}$, modulo the relations
\[
\begin{array}{rclr}
\xi_i^2 & = & 0, & i\in \N_\Box\ ; \\
\xi_\Box \xi_1 & = & \xi_1 \xi_\Box, \\
\xi_\Box \xi_2 \xi_\Box & = & \xi_2 \xi_\Box \xi_2, \\
\xi_i\xi_{i+1}\xi_i & = & \xi_{i+1}\xi_i\xi_{i+1}, & i>0\ ; \\
\xi_i\xi_j & = & \xi_j\xi_i, & j> i+1, \ \text{and} \ (i,j) \neq (\Box,2).
\end{array}
\]
As in the previous sections, for any $w\in \wt{W}_n$, choose a reduced
word $a_1\cdots a_\ell$ for $w$ and define $\xi_w := \xi_{a_1}\ldots
\xi_{a_\ell}$. The $\xi_w$ form a free $\Z$-basis of $\wt{\cW}_n$, and
we denote by $\langle \zeta,w\rangle$ the coefficient of $\xi_w$ in
the expansion of the element $\zeta\in \wt{\cW}_n$.

Let $t$ be an independent variable, define
\[
D(t) := (1+t \xi_{n-1})\cdots (1+t \xi_2)(1+t \xi_1)(1+t \xi_\Box)
(1+t \xi_2)\cdots (1+t \xi_{n-1}),
\]
and let $D(Z)=D(z_1)D(z_2)\cdots$.
Choose an integer $k$ with $0\leq
k<n$.  For any $w\in W_n$, the (restricted) type D {\em mixed Stanley
  function} $I_w(Z\,;X_k)$ is defined by
\[
I_w(Z\,;X_k) := \langle D(Z)A(X_k), w\rangle. 
\]
The power series $I_w(Z\,;X_k)$ is symmetric in the $Z$ and $X_k$
variables, separately, and has nonnegative integer coefficients.

When $w$ is the $k$-Grassmannian element associated to a typed
$k$-strict partition $\la$, then we have
\begin{equation}
\label{IE}
I_w(Z\,;X_k)=\Eta_\la(Z\,;X_k).
\end{equation}
We say that an element $w=(w_1,\ldots,w_n)$ has a {\em descent} at
position $i\in \N_\Box$ if $\ell(ws_i)<\ell(w)$.  If $k\geq 2$, we say
that $w$ is {\em increasing up to $k$} if it has no descents less than
$k$; this means that $|w_1|< w_2<\cdots < w_k$. By convention we agree
that every element of $\wt{W}_\infty$ is increasing up to $\Box$ and
also increasing up to $1$.  We can now state the following
generalization of equality (\ref{IE}): for any element $w\in \wt{W}_n$
which is increasing up to $k$, we have
\begin{equation}
\label{Ieq}
I_w(Z\,;X_k) = \sum_{\la}d^w_\la \,\Eta_\la(Z\,;X_k)
\end{equation}
summed over typed $k$-strict partitions $\la$ with $|\la|=\ell(w)$.
The integers $d^w_\la$ are nonnegative, in other words, the function
$I_w(Z\,;X_k)$ is {\em eta positive} when $w$ is increasing up to
$k$. A combinatorial interpretation for these coefficients is
provided below.

For any $w\in \wt{W}_\infty$ which is increasing up to $k$, we
construct the $k$-transition tree $\wt{T}^k(w)$ with nodes given by
elements of $\wt{W}_\infty$ and root $w$ in a manner parallel to  
Section \ref{acC}. Let $r$ be the largest descent of $w$. If $w=1$, or $k\neq 1$
and $r=k$, or $k=1$ and $r\in \{\Box,1\}$, then set $\wt{T}^k(w)
:=\{w\}$. Otherwise, let $s := \max(i>r\ |\ w_i < w_r)$ and define
$\wt{\Phi}(w) := \wt{\Phi}_1(w)\cup \wt{\Phi}_2(w)$, where
\begin{gather*}
\wt{\Phi}_1(w) := \{wt_{rs}t_{ir}\ |\ 1\leq i < r \ \ \mathrm{and} \ \ 
\ell(wt_{rs}t_{ir}) = \ell(w) \}, \\
\wt{\Phi}_2(w) := 
\{wt_{rs}\ov{t}_{ir}\ |\ i\neq r \ \ \mathrm{and} \ \
\ell(wt_{rs}\ov{t}_{ir}) = \ell(w) \}.
\end{gather*}
To define $\wt{T}^k(w)$, we join $w$ by an edge to each $v\in
\wt{\Phi}(w)$, and attach to each $v\in \wt{\Phi}(w)$ its tree
$\wt{T}^k(v)$. Then $\wt{T}^k(w)$ is a finite tree called the {\em
  $k$-transition tree} of $w$, and its leaves are all $k$-Grassmannian
elements of $\wt{W}_\infty$.  The coefficient $d^w_\la$ in (\ref{Ieq})
is equal to the number of leaves of $\wt{T}^k(w)$ of shape $\la$.

\subsection{Cohomology of flag manifolds}

Let $\{e_1,\ldots,e_{2n}\}$ denote the standard orthogonal basis of
$E:=\C^{2n}$ and let $F_i=\langle e_1,\ldots, e_i\rangle$ be the
subspace spanned by the first $i$ vectors of this basis, as in Section
\ref{cGD}.  The group $G=\SO_{2n}$ acts on the space of all complete
isotropic flags in $E$ with two orbits, determined by the family of
the maximal isotropic subspace in a given flag.  The stabilizer of the
flag $F_\bull$ is a Borel subgroup $B$ of $G$. Let $T\subset B$ denote
the maximal torus of diagonal matrices in $G$, and the Weyl group
$W=N_G(T)/T\cong \wt{W}_n$.

The parabolic subgroups $P$ of $\SO_{2n}$ with $P\supset B$
correspond to sequences $a_1<\cdots < a_p$ of elements of $\N_\Box$
with $a_p<n$. For any such $P$, the manifold $\X:=\SO_{2n}/P$
parametrizes partial flags of subspaces
\[
E_\bull \ : \ 0 \subset E_p \subset \cdots \subset E_1 \subset E=\C^{2n}
\]
with $E_1$ isotropic, $\dim(E_r) = n-a_r$ for each $r\in [1,p]$, and
$E_1$ in a given family if $a_1=\Box$. The associated parabolic
subgroup $W_P$ of $\wt{W}_n$ is generated by the simple reflections
$s_i$ for $i\notin\{a_1,\ldots, a_p\}$.  As usual, $E_r$ for $r\in
[1,p]$ and $E$ will also denote the corresponding tautological vector
bundles over $\X$.

There is a canonical presentation of the cohomology ring of $\SO_{2n}/B$,
which gives geometric significance to the variables which appear in
Section \ref{spsD}. Let $\IB^{(n)}$ denote the ideal of 
$\Gamma'[X_n]_\Q:=\Gamma'[X_n]\otimes_\Z\Q$
generated by the homogeneous elements of positive degree in
$B^{(n)}$, so that $\IB^{(n)}=\langle
{}^nb_n',{}^nb_1,{}^nb_2,\ldots\rangle$. We then have ring isomorphism
\begin{equation}
\label{BorelD}
\HH^*(\SO_{2n}/B,\Q) \cong \Gamma'[X_n]_\Q/\IB^{(n)}
\end{equation}
which maps each variable $b_p$ to $c_p(E/E_n)/2$ and $x_i$ to
$c_1(E_{n+1-i}/E_{n-i})$ for $1\leq i \leq n$. Furthermore, for any
parabolic subgroup $P$ of $\SO_{2n}$, the projection map $\SO_{2n}/B\to
\SO_{2n}/P$ induces an injection $\HH^*(\SO_{2n}/P)\hookrightarrow
\HH^*(\SO_{2n}/B)$ on cohomology rings, and we have
\begin{equation}
\label{BorelPD}
\HH^*(\X,\Q) \cong \Gamma'[X_n]_\Q^{W_P}/\IB^{(n)}_P,
\end{equation}
where $\Gamma'[X_n]_\Q^{W_P}$ denotes the $W_P$-invariant subring of
$\Gamma'[X_n]_\Q$, and $\IB_P^{(n)}$ is the ideal of $\Gamma'[X_n]_\Q^{W_P}$ 
generated by ${}^nb_n',{}^nb_1,{}^nb_2,\ldots$

We have a decomposition
\[
\SO_{2n} = \bigcup_{w \in W^P}Bw P
\]
where 
\[
W^P := \{w\in \wt{W}_n\ |\ \ell(ws_i) = \ell(w)+1,\  \forall\, i \notin
\{a_1,\ldots, a_p\}, \ i<n\}
\]
is the set of minimal length $W_P$-coset representatives in $\wt{W}_n$.
For each $w \in W^P$, the $B$-orbit of $wP$ in $\SO_{2n}/P$ is the
{\em Schubert cell} $\Y^\circ_w:=Bw P/P$. The {\em Schubert variety}
$\Y_w$ is the closure of $\Y^\circ_w$ in $\X$. Then $\X_w:=\Y_{w_0w}$
has codimension $\ell(w)$ in $\SO_{2n}/P$, and its cohomology class 
$[\X_w]$ is a {\em Schubert class}. We have an isomorphism of
abelian groups
\[
\HH^*(\X,\Z) \cong \bigoplus_{w\in W^P}\Z[\X_w]
\]
which generalizes (\ref{OGrasseq}). 

Recall that $E_r$ for $r\in [1,p]$ and $E$ denote the tautological and
trivial vector bundles over $\X$, of rank $n-a_r$ and $2n$,
respectively. For any $w\in W^P$, we then have
\begin{equation}
\label{Dpartfl}
[\X_w] = \sum_{\underline{\la}} g^w_{\underline{\la}}\,
\Eta^{(a_1)}_{\la^1}(E-E_1) s_{\la^2}(E_1-E_2)\cdots s_{\la^p}(E_{p-1}-E_p)
\end{equation}
in $\HH^*(\X,\Z)$, where the sum is over all sequences of
partitions $\underline{\la}=(\la^1,\ldots,\la^p)$ with $\la^1$ being 
typed $a_1$-strict, and the 
coefficients $g^w_{\underline{\la}}$ are given by 
\begin{equation}
\label{Dcoeff}
g^w_{\underline{\la}} := \sum_{u_1\cdots u_p = w}
d_{\la^1}^{u_1}c_{\la^2}^{u_2}\cdots c_{\la^p}^{u_p}
\end{equation}
summed over all factorizations $u_1\cdots u_p= w$ such that
$\ell(u_1)+\cdots + \ell(u_p)=\ell(w)$, $u_j\in S_\infty$ for $j\geq
2$, and $u_j(i)=i$ for all $j>1$ and $i \leq a_{j-1}$.  The
nonnegative integers $d^{u_1}_{\la^1}$ and $c^{u_i}_{\la^i}$ which
appear in the summands in (\ref{Dcoeff}) are the same as the ones in
equations (\ref{Ieq}) and (\ref{Geq}), respectively.  When $p=1$, the
partial flag manifold $\X$ is the orthogonal Grassmannian
$\OG(n-a_1,2n)$, and formula (\ref{Dpartfl}) specializes to equation
(\ref{GiambD}). In general, the polynomial
$\Eta^{(a_1)}_{\la^1}(E-E_1)$ in (\ref{Dpartfl}) is defined by pulling
back the polynomial $\Eta^{(a_1)}_{\la^1}(c(E/E_1))$ under the natural
projection map $E_\bull \mapsto E_1$ from $\X$ to $\OG(n-a_1,2n)$.

\begin{example}
\label{DSex}
Let $P=B$ be the Borel subgroup, so that the flag manifold
$\SO_{2n}/B$ parametrizes flags of subspaces $0\subset E_n \subset
\cdots \subset E_1 \subset E=\C^{2n}$ with $E_1$ isotropic in a given
family, and $\dim(E_i) = n+1-i$ for each $i\in [1,n]$. For any strict
partition $\la$, define the {\em $\wt{P}$-polynomial}
$\wt{P}_\la(X_n):=2^{-\ell(\la)}\,\wt{Q}_\la(X_n)$. Since $a_1=0$, we
have $\Eta^{(a_1)}_{\la}(E-E_1) = \wt{P}_\la(X_n)$, where we set
$x_i:=-c_1(E_i/E_{i+1})$ for each $i$. For any element $w\in
\wt{W}_n$, define the {\em orthogonal Schubert polynomial}
$\DS_w(X_n)$ by
\[
\DS_w(X_n):= \sum_{v,\om,\la}d^v_\la\,\wt{P}_\la(X_n)\AS_\om(-X_n)
\]
where the sum is over all factorizations $v\om = w$ and strict
partitions $\la$ such that $\ell(v)+\ell(\om)=\ell(w)$, $\om\in S_n$,
and $|\la|=\ell(v)$. Working as in Example \ref{CSex}, one can show
that formula (\ref{Dpartfl}) is equivalent to the statement that for
any element $w\in \wt{W}_n$, we have $[\X_w]=\DS_w(X_n)$ in
$\HH^*(\SO_{2n}/B,\Z)$.
\end{example}

\section{Historical notes and references}
\label{notes}

\subsection{Schur polynomials}

The Schur polynomials $s_\la(X_n)$ were first defined in the early
19th century by Cauchy \cite{C} using formula (\ref{Cauchydef}), as a
quotient of two alternant determinants. These polynomials were studied
further by Jacobi \cite{J} and his student Trudi \cite{Tr}, who
established the determinantal formula (\ref{equality1}). The dual
identity (\ref{JTdual}) was proved by N\"agelsbach \cite{N}. The
reformulation of the Jacobi-Trudi identity using raising operators
(\ref{giambelliA}) originates in Young's work on the representation
theory of the symmetric group; see \cite[Eqn.\ (I)]{Y} and compare
with \cite[Eqn.\ (2.26)]{R}. 

The Pieri and Giambelli formulas (\ref{PieriA}) and (\ref{GiambA}) for
Grassmannians were proved in \cite{Pi} and \cite{G}, respectively.
In his thesis \cite{S1}, Schur showed that the polynomials
$s_\la(X_n)$ can be viewed as the characters of the irreducible
polynomial representations of the general linear group $\GL_n$. From
the perspective of representation theory, equation (\ref{Cauchydef})
is a special case of the Weyl character formula. Finally, the tableau
formula (\ref{tabform}) for Schur polynomials was established by
Littlewood \cite{Li1}, more than 100 years after Cauchy's definition
appeared. For an approach to the theory of Schur polynomials starting
from the raising operator definition (\ref{giambelliA}), which
includes the above topics and more, see \cite{T4}.

The code $\gamma(\om)$ of a permutation $\om$ is the Lehmer code used
in computer science, which was known in the 19th century (see for
example \cite{La}). The shape $\la(\om)$ of a permutation was defined
in \cite{LS1}. The Stanley symmetric functions $G_\om(X)$ were
introduced in \cite{St}; in the notation of op.\ cit., the function
$G_{\om^{-1}}(X)$ is assigned to $\om$. Stanley's conjecture that the
coefficients $c^{\om}_\la$ in (\ref{Geq}) are nonnegative was proved
independently by Edelman and Greene \cite{EG} and Lascoux and
Sch\"utzenberger \cite{LS2}; the latter authors introduced the
transition trees of Section \ref{ac}. The definition of $G_\om(X)$
using the nilCoxeter algebra is due to Fomin and Stanley \cite{FS}.

The presentations (\ref{BorelA}) and (\ref{BorelPA}) of the cohomology
of type A flag manifolds are due to Borel \cite{Bo}. Formula
(\ref{Apartfl}) was proved in \cite{BKTY}, and grew out of a study of
the Schubert polynomials $\AS_\om(X_n)$ of Lascoux and Sch\"utzenberger
\cite{LS1, M1} and their relation to the quiver polynomials of Buch
and Fulton \cite{BF}. The definition of $\AS_\om(X_n)$ given in
formula (\ref{FSeq}) is found in \cite{FS}.

\subsection{Theta polynomials}

In order to clarify the relevant history, the arXiv announcement years
are listed for the main papers below, since their publication dates in
journals have little to do with when the work was completed. The story
begins with the companion papers \cite{BKT1} (arXiv:2008) and
\cite{BKT2} (arXiv:2008) which studied Schubert calculus on non
maximal isotropic Grassmannians.  The first paper \cite{BKT1} proved
the Pieri rule (\ref{PieriC}), while \cite{BKT2} dealt with the
Giambelli formula (\ref{GiambC}) and theta polynomials.

The parametrization of the Schubert classes on $\IG(n-k,2n)$ by
$k$-strict partitions was introduced in \cite{BKT1}. As explained in
op.\ cit., although the $k$-strict partitions are not really needed
there, they are a key ingredient of its companion paper \cite{BKT2},
and in related works such as \cite{T3} (arXiv:2008). 

The paper \cite{BKT2} was the first to realize that Young's raising
operators play an essential role in geometry, in the Giambelli type
formulas for isotropic Grassmannians, and to employ them in their
proofs. Before \cite{BKT2}, these operators made occasional
appearances, notably in the theory of Hall-Littlewood polynomials (see
e.g.\ \cite{Li2, Mo, M2}), but they were rarely used, even in
representation theory and combinatorics. The solution to the Giambelli
problem for the usual (type A), Lagrangian, and maximal orthogonal
Grassmannians found in \cite{G} and \cite{P}, respectively, employed
the older language of Jacobi-Trudi determinants (in Lie type A) and
Schur Pfaffians (in types B, C, and D), which goes back to \cite{J} and
\cite{S2}.

In the initial version of \cite{BKT2}, the theta polynomials were
expressed as the formal power series in (\ref{Thetatheta}), although
the intention in op.\ cit.\ (which justifies the term `polynomial')
was to regard the $\Ti_\la$ as Giambelli polynomials in the
$\ti_p$. The definition (\ref{giambelliC}) of the theta polynomials
$\Ti_\la(u)$ in independent variables $u_i$ was first given explicitly
in \cite[Eqn.\ (3)]{T5} (arXiv:2009).

The ring $\Gamma$ of Section \ref{spsC} is isomorphic to the ring of
Schur $Q$-functions (see Example \ref{exQ}, \cite{S2}, and
\cite[III.8]{M2}), whose elements are symmetric formal power series in
$Z:=(z_1,z_2,\ldots)$. Using the latter notation, the action of
$W_\infty$ on the ring $\Gamma[X]$ was studied by Billey and Haiman in
\cite[Lemma 4.4]{BH}.  The same authors obtained a natural $\Z$-basis
of $\Gamma[X]$ consisting of type C Schubert polynomials (actually
power series) $C_w(X)$, for $w\in W_\infty$. When $w$ is an
$n$-Grassmannian element of $W_\infty$ of shape given by the
$n$-strict partition $\la$, then $C_w(X)$ is equal to $\Ti_\la(X_n)$
in $\Gamma[X]$, so modulo the defining relations in $\Gamma$ (see
\cite[Prop.\ 6.2]{BKT2}).

It was shown in \cite{T9} that the theta polynomials $\Ti_\la(X_n)$
of level $n$ are symmetric for the $W_n$-action on $\Gamma[X_n]$, and
form a $\Z$-basis for the $W_n$-invariants there, which is the content
of (\ref{Cinv}) and (\ref{GTieq}). The paper \cite{T9} goes on to
define the shape $\la(w)$ of a signed permutation $w$, proves equation
(\ref{mScheq}), and also obtains the presentations (\ref{BorelC}) and
(\ref{BorelPC}).  We note that multi-Schur Pfaffians such as
(\ref{multiSPf}) first appeared in the work of Kazarian \cite{K},
resurfaced in \cite{IMN} (arXiv:2008), and were subsequently used in
the degeneracy locus formulas of \cite{AF}.

The tableau formula (\ref{tableauxeq}) for $\Ti_\la(Z\,;X_k)$ was
established in \cite{T3}. The fact that the {\em left weak Bruhat
  order} on the $k$-Grassmannian elements of $W_\infty$ respects the
inclusion relation $\la\subset\mu$ of $k$-strict Young diagrams is
what makes such a formula possible. This was pointed out in
\cite[Prop.\ 4]{T3}. Example \ref{Qex} is well known in the 
combinatorial theory of Schur $Q$-functions; see for instance 
\cite[III.(8.16)]{M2}.

The type C mixed Stanley functions and $k$-transition trees were
defined and equation (\ref{Jeq}) was proved in \cite{T5}. These
constructions build on earlier results of Billey and Haiman \cite{BH,
  B} and Fomin and Kirillov \cite{FK}, which had studied the $k=0$
case. The paper \cite{T5} also proves formula (\ref{Cpartfl}), which
gives an intrinsic solution to the Giambelli problem of representing
the Schubert classes in the cohomology of $\Sp_{2n}/P$, for any
parabolic subgroup $P$ of $\Sp_{2n}$. 

The $\wt{Q}$-polynomials $\wt{Q}_\la(X_n)$ in Example \ref{CSex} are
due to Pragacz and Ratasjki \cite{PR}. The symplectic Schubert
polynomials $\CS_w(X_n)$ in (\ref{Teq}) were defined in \cite{T1}
(arXiv:2008), and are a {\em geometrization} of the Billey-Haiman type
C Schubert polynomials. For more on this history, we refer the reader
to \cite[Section 5]{T8}.

\subsection{Eta polynomials}

The paper \cite{BKT1} introduced typed $k$-strict partitions and
proved the Pieri rule (\ref{pruleD2}) for non maximal even orthogonal
Grassmannians. Young's raising operators were used in \cite{BKT4}
(arXiv:2011) to define eta polynomials and to prove equation
(\ref{GiambD}), which solves the Giambelli problem for the same
spaces. The work \cite{T5} extended this to address the analogous
question for all partial orthogonal flag manifolds, which gives
equation (\ref{Dpartfl}). See \cite{T7} for a detailed exposition,
which includes the version of this result which holds in the more
general setting of degeneracy loci of vector bundles for the classical
Lie groups.

The ring $\Gamma'$ of Section \ref{spsD} is isomorphic to the ring
$\Z[P_1,P_2,\ldots]$ of Schur $P$-functions. The action of
$\wt{W}_\infty$ on the ring $\Gamma'[X]$ and the induced divided
difference operators there were studied by Billey and Haiman
\cite{BH}, who defined a $\Z$-basis of $\Gamma'[X]$ consisting of type
D Schubert polynomials $D_w(X)$ for $w\in \wt{W}_\infty$, compatible
with these operators.  According to \cite[Prop.\ 6.3]{BKT4}, when $w$
is an $n$-Grassmannian element of $\wt{W}_\infty$ associated to the
typed $n$-strict partition $\la$, then $D_w(X)$ is equal to
$\Eta_\la(X_n)$ in $\Gamma'[X]$. Again it is important to note that
this equality takes place in a ring with relations, which come from
the subring $\Gamma'$.

The fact that the eta polynomials of level $n$ are symmetric for the
$\wt{W}_n$-action on $\Gamma'[X_n]$ and provide a $\Z$-basis for the
Weyl group invariants there, which is the content of (\ref{Binv}) and
(\ref{GTieqD}), was explained in \cite{T9}.  The paper \cite{T9} also
defines the shape $\la(w)$ of an element $w$ of $\wt{W}_\infty$,
proves equation (\ref{mScheqD}), and moreover obtains the
presentations (\ref{BorelD}) and (\ref{BorelPD}).

The tableau formula (\ref{tableauxeqD}) for the eta polynomials
$\Eta_\la(Z\,;X_k)$ was established in \cite{T6} (arXiv:2011). The
type D mixed Stanley functions and $k$-transition trees were defined
and equation (\ref{Ieq}) was proved in \cite{T5}.  This used the even
orthogonal (type D) version of the nilCoxeter algebra approach of
\cite{FS, FK} to type A and B Stanley symmetric functions, which is
found in \cite{Lam}. The $\wt{P}$-polynomials $\wt{P}_\la(X_n)$ and
orthogonal Schubert polynomials $\DS_w(X_n)$ of Example \ref{DSex}
were defined in \cite{PR} and \cite{T2} (arXiv:2009), respectively.

\end{document}